\documentclass[12pt]{article}
\usepackage{amsfonts,supertabular,latexsym,amsmath}

\newtheorem{thm}{Theorem}[subsection]

\newcommand{\F}{\mbox{$\mathbb F$}}
\newcommand{\Q}{\mbox{$\mathbb Q$}}     
\newcommand{\R}{\mbox{$\mathbb R$}}     
\newcommand{\C}{\mbox{$\mathbb C$}}     
\newcommand{\Z}{\mbox{$\mathbb Z$}}     

\renewcommand{\exp}{{\rm exp}\,}

\renewcommand{\tilde}{\widetilde}

\renewcommand{\bar}{\overline}
\newcommand{\shom}{\mbox{${\cal H}om\,$}}  

\newcommand{\veq}{\mbox{\large $\parallel$}}  
\newcommand{\sZ}{\mbox{\scriptsize{$\Z$}}}    
\newcommand{\sF}{\mbox{\scriptsize{$\F$}}}   

\newcommand{\boxtensor}{{\Box\kern-9.03pt\raise1.42pt\hbox{$\times$}}}

\newcommand{\arc}[1]{{#1}\kern-16pt\raise7pt\hbox{\Large{$\frown$}}}

\newcommand{\longinto}{\lhook\joinrel\kern-3pt\hbox to
100pt{\rightarrowfill}}           

\newcommand{\propsubset}
{\mbox{$\textstyle{
\subseteq_{\kern-5pt\raise-1pt\hbox{\mbox{\tiny{$/$}}}}}$}}

\newcounter{elno}                

\newcounter{example}[section]
\def\theexample{\thesection.\arabic{example}}

\newcounter{exercise}[section]
\def\theexercise{\thesection.\arabic{exercise}}

\def \psb{parabolic stable bundle}
\def \pssb{parabolic semistable bundle}

\def \phn{parabolic Harder Narasimhan}
\def \pa{parabolic}
\def \ex{extension}
\def \la{\lambda} 
\def \ps={`par semi-stable = par stable'}
\def \vsp{{\vskip .5cm}}
\def \vs{{\vskip 1.5mm}}
\title{Poincar\'e polynomial of the moduli spaces of parabolic bundles}
\author{Yogish I. Holla}
\begin{document}
\date{Jan 15 1999}
\maketitle
School of Mathematics, Tata Institute of Fundamental Research,
Homi \\ Bhabha Road, Mumbai 400 005, India. 
e-mail: yogi@math.tifr.res.in

\begin{abstract}In this paper we use Weil conjectures (Deligne's
theorem) to calculate the Betti numbers of the moduli spaces of
semi-stable parabolic bundles on a curve. The quasi parabolic analogue
of the Siegel formula, together with the method of Harder-Narasimhan
filtration gives us a recursive formula for the Poincar\'e
polynomials of the moduli. We solve the recursive formula by the
method of Zagier, to give the Poincar\'e polynomial in a closed form.
 We also give explicit tables of Betti numbers in small rank, and genera.
\end{abstract}

\section{Introduction}
 
This paper uses the Riemann hypothesis of Weil (Deligne's theorem) to 
explicitly determine the Betti numbers of the moduli of semistable 
parabolic bundles on a curve (when parabolic semi-stability implies
parabolic stability).
\vs
Vector bundles with parabolic structures were introduced by Seshadri,
and their moduli was constructed by Mehta-Seshadri (see [S] for an account). 
Our approach to the calculation of the Betti numbers is an extension
of the method used by Harder and Narasimhan [H-N] in the case of ordinary
vector bundles. Harder and Narasimhan use the result of Siegel that 
the Tamagawa number of $SL_r$ over a function field of trancedence 
degree one over a finite field is $1$. This result can be reformulated
purely in terms of vector bundles to give the formula (equation (2.16)),
which was used by Desale and Ramanan [D-R] in their refinement of the 
Harder-Narasimhan Betti number calculation. 
\vs
In place of the above formula, we use its quasi-parabolic analogue 
(see equation (2.19)) proved by Nitsure [N2], to extend the
calculation of Harder and Narasimhan, as refined by Desale and Ramanan,
 to parabolic case. 
\vs    
This gives us a recursive formula to obtain Betti numbers. Such a 
recursive formula had been obtained earlier for genus $\ge$ 2 by
Nitsure[N1] using the Yang-Mills method of Atiyah-Bott[A-B], 
and this was extended to lower genus by Furuta and Steer[F-S]. 
\vs
Finally following Zagier's[Z] method of solving such
a recursion(in the case of ordinary vector bundles), we obtain an
explicit formula for the Poincar\'e polynomials. We give sample tables
in lower ranks and genera.
\vs
This paper is arranged as follows. 
In section 2, we have introduced our notations and recalled certain
basic facts about parabolic bundles for the convenience of the reader.
The paper of Desale and Ramanan computes the Poincar\'e polynomial
of the moduli space of stable bundles, starting with the formula of
Siegel (2.16). In section 3, we have followed their
general pattern with suitable changes needed to handle the parabolic
case, with the Siegel formula replaced by its parabolic analogue
(2.19). This gives us the theorem (3.36), which is our desired
recursive formula for the Poincar\'e polynomial. Along the way, we need a
certain substitution ($\omega_i \rightarrow -t^{-1},~q \rightarrow
t^{-2}$) used by Harder and Narasimhan, who have sketched its
justification. We give a detailed proof of why such a substitution
works (in a somewhat more general context) in
section 4. In section 5, we solve the recursive formula using
Zagier[Z]'s approach, to get the explicit form (5.23) of the
Poincar\'e polynomial. In section 6,
we give some sample computations of the Poincar\'e polynomials of
these moduli spaces and check their dependence on the weights and the
degree when the rank is low (2,3 and 4). In the appendix (section 7), we
have given tables for the Betti numbers of these moduli spaces in 
rank 2,3 and 4. 
\vs
{\bf Acknowledgments}\, The author wishes to thank Nitin Nitsure 
for his constant support and encouragement and for many useful
comments and suggestions. The author also wishes to thank
K.P.Yogendran for writing the computer program which was used to
calculate the Betti numbers for rank 3 and 4.

\section{Basic definitions and notations}

\centerline {\bf Zeta function of a curve.}
\medskip
~~~\,~Let ${\F}_q$ be a finite field, and let $\bar{\F}_q$ be its
algebraic closure.
Let $X$ be a smooth projective geometrically irreducible curve over
${\F}_q$, where geometric irreducibility means 
$\bar{X}=X \otimes_{ \sF_q}{\bar{\F}_q}$ is irreducible.
\vs
  Given any integer $r >0$, let ${\F}_{q^r} \subset \bar{\F}_q $ be 
the unique field
extension of degree $r$ over $\F_q$. Let $N_r = |X(\F_{q^r})|$ be the 
cardinality of the set of $\F_q$-rational points of $X$. 
  Recall that the zeta function of $X$ is defined by
$$ 
Z_X(t) = \exp \left({\sum_{r > 0}\frac{ N_rt^r}{r}}\right)  \eqno (2.1)
$$
 By the Weil conjectures it follows that the zeta function has the
 form
$$  Z_X(t) =  \frac{\prod_{i=1}^{2g}(1 - \omega_it)}{(1 - t)(1 -
qt)}          \eqno(2.2)
$$             
 where $\omega_i$'s are algebraic integers of norm $q^{1/2}$, and $g$ is
the genus of the curve.
For $\nu \ge 1$, let $X_{\nu}$ denote the curve 
 $X_{\nu} = X \otimes_{ \sF_q }{ \F_{q^{ \nu}}}$.
 The following remark will be used later.
\vsp
{\bf Remark\, 2.3.}\,
 If the zeta function of $X$ over $\F_q$ is as given in
(2.2), then the zeta function for the curve $X_{\nu}$ over $\F_{q^{\nu}}$
has the form 
$$ 
Z_{X_{\nu}}(t) = \frac{\prod_{i=1}^{2g}(1 - \omega_i^{\nu}t)}{(1 - t)(1 -
qt)}.       \eqno (2.4)
$$
\vsp
\medskip
\centerline {\bf Rational points on Flag varieties}
\medskip
Now we recall the computation of the number of rational points of flag
varieties. Let $k=\F_q$ as before,
let $n$ and $m$ be positive integers, and let there
be given non-negative integers $r_1,\ldots ,r_m$ with
$r_1+\ldots +r_m =n$. We denote by Flag$(n,m,(r_j))$ the variety
of all flags
$k^n=F_1\supset\ldots\supset F_m\supset F_{m+1}=0$
of vector subspaces in $k^n$, with
dim$(F_j/F_{j+1}) = r_j$.
\vsp
{\bf Proposition\, 2.5.}\,
{\it The  number of $\F_q$-rational points of $\mbox{Flag}(n,m,(r_j))$ is
$$f(q,n,m,(r_j)) = \frac {\prod_{i=1}^n{(q^i -
1)}}{\prod_{\{j|r_j \neq 0\}}\prod_{l=1}^{r_j}{(q^l -1)}}   \eqno(2.6)$$
}
\vsp
{\noindent{\bf Proof.}} The number of rational points $g(r,p)$ on the 
Grassmanian Grass$(r,p)$ of $p$-dimensional subspaces of $k^r$ can be
seen to be
$$g(r,p) = \frac{(q^r-1)\cdots(q^r-q^{p-1})}{(q^p-1)\cdots
(q^p-q^{p-1})}    \eqno(2.7)$$
and the number $f(q,n,m,(r_j))$ clearly satisfies
$$f(q,n,m,(r_j))= g(n,r_m)g(n-r_m,r_{m-1})\cdots
g(r_1+r_2,r_2)   \eqno(2.8)$$
Simplifying yields the desired formula.   $\hfill{\Box}$
\vsp
\medskip
\centerline {\bf Parabolic vector bundles}          
\medskip
Let $S=\{ P_1,\ldots, P_s\}$ be any closed subset of $X$ whose points
are $k$-rational. For each $P \in S$, let there be given a positive 
integer $m_P$.
\vs
We fix an indexed family of real numbers $( \alpha^P_i )$, where 
$P\in S$ and $i = 1, \ldots ,m_P$, 
satisfying $0\le \alpha^P_1
< \alpha^P_2 \ldots < \alpha^P_{m_P}  < 1$,
which we denote simply by $\alpha$. 
We fix the set $S$, the integers $(m_P)$ and the family 
$\alpha$ in all that follows.
\vs
The parabolic weights at $P$
of the parabolic bundles that we will
consider in this paper are going to belong to the chosen set
$\{ \alpha^P_1,\alpha^P_2, \ldots ,\alpha^P_{m_P} \}$.
(Note that this property will be inherited by the 
sub-quotients of such parabolic bundles.) 
This allows us to formulate the definition of parabolic bundles in a 
somewhat different way from Seshadri, which is more suited for our
inductive arguments. However, the difference is only superficial, and
we explain later (remark (2.11) on the next page) the 
bijective correspondence between parabolic bundles in our sense, and
parabolic bundles in Seshadri's sense which have weights in our given set.
\vs
A {\bf quasi-parabolic data} $R$ (or simply `data' when the context is
 clear) is an indexed family of non-negative integers $(R^P_i)$
for $P\in S$ and $1 \le i \le m_P $, satisfying the following condition:
$\sum_{i=1}^{m_P}R^P_i$ is a positive integer independent of $P\in S$.
We call $n(R)=\sum_{i=1}^{m_P}R^P_i$ as the {\bf rank} of the
 quasi-parabolic data $R$.
\vs
Let $L$ be another quasi-parabolic data. We say $L$ is a {\bf
sub-data} of a given data $R$
if $ L^P_i \le R^P_i ~\mbox{for all}~P~ \mbox{and}~i $, and $ n(L) < n(R)$.  
We also define its {\bf complementary sub-data} $R-L$
by $(R-L)^P_i = R^P_i -  L^P_i $. 
 \vs 

A quasi-parabolic structure with data $R$ on a vector bundle $E$ on $X$, by
definition, consists of a flag
\vs
  $$E^P = E^P_1 \supset 
~E^P_2 \ldots \supset ~E^P_{m_P} \supset ~E^P_{m_P+1}= 0 \eqno(2.9)$$ 
of vector subspaces in the fiber $E^P$  over each  point $P$ of $S$,
such that $R^P_i = \dim(E^P_i/E^P_{i+1})$ for $P \in S$ and $1 \le i \le m_P$.
 \vs


A parabolic structure with data $R$ on a vector bundle $E$  is a
quasi-parabolic structure on $E$ with data $R$ along 
with weights $\alpha^P_i$ for each $P$ and $i$. We say
$R^P_i$ is the multiplicity of the weight $\alpha^P_i$.
To a data $R$, we associate the real number $\alpha(R)$ by    
$$\alpha(R) = \sum_P \sum_{i=1}^{m_P}{R^P_i \alpha^P_i}. \eqno(2.10)$$
\vsp
{\bf Remark\, 2.11.}\,
 We record here the minor changes in notations and
conventions that we have made (compared to the original notation of Seshadri).
In our definition, note that the data $R$ has the property that $R^P_i
\ge 0$ (and not $> 0$), hence if $E$ is a parabolic bundle in our
sense with data $R$ then the inclusions occurring in the filtration 
(2.9) are not necessarily strict. We recover the 
definition of  Seshadri by re-defining the weights 
${\bar \alpha}$ inductively as follows 
$$ 
\begin{array}{ll}
{\bar \alpha}^P_1 &=~~ \underset {i}{\mbox{min}}
                         ~\{\alpha ^P_i| R^P_i \neq 0\} \\
{\bar \alpha}^P_j &=~~ \underset {k}{\mbox{min}}~\{\alpha ^P_k| R^P_k \neq 0,
\alpha^P_k - {\bar \alpha}^P_{j-1} > 0 \}. 
\end{array}
\eqno(2.12)
$$
From this it follows that each ${\bar \alpha}^P_j$ equals $\alpha ^P_i$ for
exactly one $i$, which allows us to define ${\bar R}^P_j= R^P_i$
for that particular $i$.
 Now it is clear that $E$ is a parabolic bundle in the sense of
Seshadri with weights ${\bar \alpha}$ and multiplicities ${\bar R}$,
where the flags are defined by sub-spaces ${\bar E}^P_j = E^P_i$ for
that $i$ for which ${\bar \alpha}^P_j = \alpha ^P_i$. 
Since we have fixed the weights $\alpha$, we can recover the parabolic
bundle in our sense
from a given parabolic bundle with weights ${\bar \alpha}$ and
multiplicities ${\bar R}$ in the sense of Seshadri when
the set $\{ {\bar\alpha}^P _i\}$ is a subset of $\{\alpha^P_i \}$ for
each $P$, 
by simply assigning 
$$
\begin{array}{ll}         
R^P_i &= 0 ~~\mbox{if}~~ \alpha^P_i \neq {\bar \alpha}^P_j
~~\mbox{for any}~~ j \\
  &= {\bar R}^P_j  ~~\mbox{if}~~  \alpha^P_i = {\bar \alpha}^P_j
~~\mbox{for some}~~ j
\end{array}          
\eqno(2.13)
$$
and the defining the sub-spaces occurring in the flags inductively by 
$$
\begin{array}{ll}
E^P_1 &= E^P\\       
E^P_i &= E^P_{i-1} ~~\mbox{if}~~ \alpha^P_i \neq {\bar \alpha}^P_j
~~\mbox{for any}~~ j \\
  &= {\bar E}^P_j  ~~\mbox{if}~~  \alpha^P_i = {\bar \alpha}^P_j
~~\mbox{for some}~~ j
\end{array}          
\eqno(2.14)
$$
This sets up a bijective correspondence between parabolic bundles
in our sense and in the sense of Seshadri. Also note that 
$$\sum_P \sum_i{R^P_i \alpha^P_i} =
\sum_P \sum_i{{\bar R}^P_i {\bar \alpha}^P_i}, \eqno(2.15)$$
which will enable us to write the parabolic degree in terms of our
modified definition.
The advantage of our definition is that it is easier to handle the
induced parabolic structures on the sub-bundles and the quotient
bundles in what follows.
Also the parabolic homomorphisms between two parabolic bundles $E$
and $E'$ with data $R$ and $R'$ respectively in our sense, having the
same fixed family of weights $(\alpha^P_i)$, are just filtration
preserving homomorphisms of the vector bundles.
\vsp
\medskip
\centerline{\bf Quasi-parabolic Siegel formula}
\medskip
For a positive
integer $n$ and for any line bundle ${\cal L}$ on $X$, let $J_n({\cal L})$ 
denote the set of isomorphism classes of vector bundles 
$E$ on $X$ with $rank(E)=n$ and determinant ${\cal L}$.
Let $|\mbox{Aut}(E)|$ denote the cardinality of the group of all
automorphisms of $E$. Then the Siegel formula, asserts that
$$\sum _{E\in J_n({\cal L})} \frac{1}{{|\mbox{Aut}(E)|}}=
\frac{q^{(n^2-1)(g-1)}}{q-1} Z_X(q^{-2})\cdots Z_X(q^{-n}) \eqno(2.16)$$
The above formula was given a proof purely in terms of vector bundles
by Ghione and Letizia [G-L]. 

For a line bundle ${\cal L}$ on $X$, let $J_R({\cal L})$ denote the set of all 
isomorphism classes of quasi-parabolic vector bundles with data $R$, 
and determinant ${\cal L}$.
Let $f_R(q)$ (denoted by $f(q,R)$ in [N2]) be the number of
   $\F_q$-valued points of the variety 
${\cal F}_R = \prod_{P \in S} \mbox{Flag}(n(R),m_P,(R^P_i))$
where $ \mbox{Flag}(n(R),m_P,(R^P_i))$ is the flag variety
   determined by $(R^P_i)$.
 Now by equation (2.6), we have
 $$f_R(q)= \frac {\prod_{i=1}^{n(R)}{(q^i -1)}^{|S|}}
{\prod_{P \in S}\prod_{\{i| {R^P_i} \ne 0\}}\prod_l^{R^P_i}{(q^l
-1)}}.   \eqno(2.17)$$
Let $|\mbox{ParAut(E)}|$ denote the cardinality of the set of quasi-parabolic
isomorphisms of a quasi-parabolic bundle $E$.
The Siegel formula has the following quasi-parabolic analogue,
which was proved by Nitsure [N2].
\vsp
{\bf Theorem\, 2.18.}\,
{\it (Quasi-parabolic Siegel formula)
$$\sum_{E\in J_R({\cal L})} \frac {1}{|\rm{ParAut}(E)|}=
f_R(q) \frac {q^{(n(R)^2-1)(g-1)}}{ q-1 } 
Z_X(q^{-2})\ldots Z_X(q^{-n(R)}) \eqno(2.19)$$
}
\vsp
For example, if $S$ is empty or more generally if the 
quasi-parabolic structure at each point of $S$ is trivial (that
is, each flag consists only of the zero subspace and the whole space), 
then on one hand $\mbox{ParAut}(E)=\mbox{Aut}(E)$, and on the other hand each flag 
variety is a point, and so $f_R(q)=1$. Hence in this situation 
 the above formula reduces to the original Siegel formula.
\vs
\medskip
\centerline {\bf Parabolic degree and stability}          
\medskip
Let $E$ be a parabolic bundle over $X$ with data $R$. 
Because of (2.15), we can define the parabolic degree of $E$ and the
parabolic slope of $E$ as follows:
  $$\mbox{pardeg}(E) = \mbox{deg} (E) + \alpha (R)~~~\mbox{and}~~~  
\mbox{par}\mu(E) = \mbox{pardeg} (E)/\mbox{rank} (E). \eqno(2.20) $$
 
A parabolic bundle $E$ on $X$ is said to be parabolic
stable (resp. parabolic semi-stable) if for every non-trivial proper
sub-bundle $F$ of $E$ with induced \\
parabolic structure,
 we have $\mbox{par}\mu(F) < \mbox{par}\mu(E) (\mbox{resp.} \le).$
 \vs
The equation (2.15) implies that the definitions of parabolic
stable(resp. parabolic semi-stable) bundles are not altered by the
change in the definition of the parabolic bundles we have made. 

We say that the numerical data $(d,R)$ satisfies the condition 
`par semi-stable = par stable' if every
 \pa ~semistable bundle with data $R$ and degree $d$ is automatically 
\pa ~stable.
\vsp
{\bf Remark\, 2.21.}\,
If the degree $d$ and rank $n(R)$ are coprime and all the weights are
 assumed to be very small ($\alpha^P_i<1/(n(R)^2|S|)$ for example) then,
 each \pssb ~ is  actually parabolic stable.
\vsp
We now recall the following.
\vsp
{\bf Lemma\, 2.22.}\,
{ \it If $E$ is a \psb ~, then every \pa ~homomorphism of
$E$ into itself is a scalar endomorphism.}
\vsp
\medskip
\centerline {\bf Parabolic Harder-Narasimhan-Intersection types}          
\medskip
Recall the following.
\vsp
{\bf Proposition\, 2.23.}\, 
{\it Any non-zero parabolic bundle $E$ with the data $R$
admits a  unique filtration by sub-bundles
$$0 = G_0 \propsubset G_1 \propsubset \ldots \propsubset G_r = E \eqno(2.24)$$
satisfying \\
{\bf  i)} $G_i/G_{i-1}$ is parabolic semi-stable for $i=1,\ldots,r$ \\
{\bf ii)}  ${\rm par}\mu(G_i/G_{i-1})>{\rm par}\mu(G_{i+1}/G_i)$ for $i=1,\ldots,r-1.$\\
   Equivalently, \\
{\bf 1)} $G_i/G_{i-1}$ is parabolic semi-stable for $i=1,\ldots,r$   \\
{\bf 2)} For any parabolic sub-bundle $F$ of $E$ containing
         $G_{i-1}$ we have \\
         ${\rm par}\mu(G_i/G_{i-1})>{\rm par}\mu(F/G_{i-1}), ~ i=1,\ldots,r.$
}
 \vsp
 The result is first proved over an algebraically closed field, and
then Galois descent is applied (using the uniqueness of the
filtration) to prove that the filtration is defined  over the original field.
\vsp
The unique filtration is called the parabolic Harder-Narasimhan
filtration. 
 \vs 
Let $(I^P_{i,k})$ be an indexed collection of non-negative integers,
where $P \in S$, $1 \le i \le m_P$, and $1 \le k \le r$ where $r$ is a
given positive integer. We say that
$I =(I^P_{i,k})$ is a {\bf partition}  of
$R$ of length $r$ if the following holds:\\
{\bf 1)} For $P\in S$ and $ 1 \le i \le m_P$, we have
$\sum_{k= 1}^r{I^P_{i,k}} = R^P_i$.\\
{\bf 2)} For $P\in S$ and $1 \le k \le r$, the summation $\sum_{i=0}^{m_P}
{ I^P_{i,k}}$
is independent of $P$.\footnote{For later reference, 
this number will be equal to the rank of $G_k$.} \\
{\bf 3)} Given any $k \le r$, $ I^P_{i,k} \ne 0$ for some $P$ and $i$.

We write $\ell (I)=r$ to indicate that $I$ has length $r$. 
\vs
Suppose $I$ is a partition of $R$ with $\ell (I)=r$. For $j = 1 ,\ldots,r$
 define a sub-data $R^I_j$  of $R$ by the
 equality 
$ (R^I_{j})^P_i = I^P_{i,j}$. 
We also define the sub-data
$R^I_{\le j}$ (resp. $R^I_{\ge j}$) of $R$ by the
 equality 
$$ (R^I_{\le j})^P_i = \sum_{k \le
j}{I^P_{i,k}}~~~(\mbox{resp.}~~~(R^I_{\ge j})^P_i= \sum_{k \ge
j}{I^P_{i,k}}).    \eqno(2.25)$$ 
Note that the rank $n(R^I_{\le j})$ of $R^I_{\le j}$ is
equal to $n(R^I_1)+n(R^I_2) \ldots + n(R^I_j)$.

We observe that the partition $I$ of $R$ induces a partition $I_{\le
j}$ (resp. $I_{\ge j}$) on $R^I_{\le j}$ (resp. $R^I_{\ge j}$) defined
by $(I_{\le j})^P_{i,k}=I^P_{i,k}$ ( resp. $(I_{\ge
j})^P_{i,k}=I^P_{i,k}$), where $k \le j$ (resp. $k \ge j$).
           
\smallskip 
We now recall how partitions, as abstractly defined above, are
associated with parabolic bundles in Nitsure[N1].
To each $E \in J_R({\cal L})$ we have the parabolic Harder-Narasimhan
filtration $ 0 \subset G_1 \subset \ldots \subset G_r = E$
 which gives a filtration on the fibers. Then the {\bf intersection
 matrix} $(I^P_{i,k})$ corresponding to it is defined in [N1], by putting  
              
$$ 
I^P_{m_P,1} = \dim(E^P_{m_P} \cap  G^P_1)   \eqno(2.26)
$$ and
             
$$ 
I^P_{j,l} = \dim( E^P_j \cap G^P_l) - \sum_{\stackrel{i \le
j,~\mbox{{\footnotesize and}}~ k \le l}{{\scriptstyle (i,k) \ne (l,j)}}}{I^P_{i,k}}  \eqno(2.27)
$$
with this definition $I$ becomes a partition of $R$ with $\ell (I)=r$.
The sub-bundles $G_j$ (resp. quotients $E/G_j$) 
under the induced parabolic structure have the
sub-data $R^I_{\le j}$ (resp. $R^I_{\ge j+1}$).  
Also the sub-quotient $G_j/G_{j-1}$ has the sub-data $R^I_j$.

\medskip
\centerline {\bf Moduli spaces}          
\medskip
For the moment assume that our ground field $k$ is algebraically closed.\\ 
Recall that \pssb s, with a fixed parabolic slope, 
form an abelian category with the property that each object has finite
length and simple objects are precisely the parabolic stable bundles. 
Hence for every parabolic semi-stable bundle $E$ there exists a 
Jordan-Holder series 
 $$E = E_r \supset E_{r-1} \ldots \supset E_1 \supset 0$$
such that $E_i/E_{i-1}$ is a \psb~ satisfying
 $\mbox{par}\mu(E_i/E_{i-1}) = \mbox{par}\mu(E)$. If we write $Gr(E)$ for
$\oplus_i{E_i/E_{i-1}}$, then it is well defined and is a \pssb ~with
the same data as $E$. We say that two \pssb s $E$ and $F$ 
are S-equivalent if $Gr(E)$ and $ Gr(F)$ are
isomorphic as parabolic bundles. 
\vs

   Mehta and Seshadri[M-S] prove that there exists a coarse moduli
scheme ${\cal M}_{R,{\cal L}}$ of the S-equivalence classes of \pssb s
with the data $R$ and determinant ${\cal L}$.
The scheme ${\cal M}_{R,{\cal L}}$ is a normal projective variety.
 Further the subset ${\cal M}_{R,{\cal L}}^s$ of ${\cal M}_{R,{\cal
L}}$ corresponding to \psb s is a smooth open subvariety.
\vsp
{\bf Remark\, 2.28.}\,
Our method computes the Betti numbers of the moduli space of parabolic
bundles for any curve over $\C$ because of the following reason. The theorem
of Seshadri implies that the topological type of these moduli spaces
depend only on the genus $g$, the cardinality of parabolic vertices $|S|$, the
degree $d$ and the set of weights $\alpha$ along with their
multiplicities $R$. We start with such a data, construct a smooth 
projective absolutely irreducible curve $X$ over a finite field 
$k=\F_q$ which has at least $|S|$ number of $k$-rational points 
(by taking $q=p^n$ for large $n$, or $q=p$ for a large prime $p$).
The use of Witt vectors allows us to spread the curve and the moduli
spaces to the quotient field of the ring of Witt vectors, when the
condition \ps= ~holds. 
Now by Weil conjectures it follows that the Betti numbers of the 
moduli space of parabolic bundles over $X$ coincides with the one 
over the curve obtained by the change of base to $\C$.

\medskip
\centerline {\bf Parabolic extensions}
\medskip
Let $E', E$ and $E''$ be parabolic bundles with data $R', R$ and
$R''$ respectively. Let $$0 {\longrightarrow} E'
\stackrel{i}{\longrightarrow} E \stackrel{j}{\longrightarrow} E''
\longrightarrow 0       \eqno(2.29)$$ 
be a short exact sequence of the underlying
vector bundles such that the parabolic structures induced on $E'$ and
$E''$ from the given parabolic structure on $E$ coincide with the
given parabolic structures of $E'$ and
$E''$, we say that (2.29) is a short exact sequence of \pa
~bundles. We also say $[E] = (E,i,j)$ is a {\bf parabolic extension} 
of $E''$ by $E'$.
\vs
 
We say two \pa ~\ex s $[E_1]$ and $[E_2]$ are equivalent if there exists
an isomorphism of \pa ~bundles $ \gamma :E_1\longrightarrow E_2$ such
that the following diagram with commutes:
     \[
     \begin{array}{lrclllllllr}
 & {0} & {\longrightarrow} & E' & \stackrel {i_1} {\longrightarrow} & E_1 &
 \stackrel {j_1} {\longrightarrow} & E'' & {\longrightarrow} & {0} &  \\
 ~~~~~~~~~~~~~ & & & \veq & & {\downarrow \gamma}  & & 
\veq & & { }    &~~~~~~~~~~~~~~~~~    (2.30)\\
 & {0} & {\longrightarrow} & E' & \stackrel {i_2} {\longrightarrow} & E_2 &
 \stackrel {j_2} {\longrightarrow} & E''  & {\longrightarrow} & {0} &    
\end{array}
\] 
 We denote the set of equivalence classes of parabolic extensions 
by \\$\mbox{ParExt}(E'',E')$. 
The proof of the following lemma is straight-forward and we omit it.
\vsp
{\bf Lemma \, 2.31.}\,
 {\it There is a canonical bijection between ${\rm ParExt}(E'',E')$
and \\$H^1(X,{\mbox{${\cal P}ar{\cal H}om\,$}}(E'',E'))$, where 
${\mbox{${\cal P}ar{\cal H}om\,$}}(E'',E'))$ is the sheaf of germs of
parabolic homomorphisms from $E''$ to $E'$.
}
\vsp
\medskip
By analogy with the case of ordinary vector bundles, we define an action of \\
$\mbox{ParAut}(E'') \times \mbox{ParAut}(E')$ on
$\mbox{ParExt}(E'',E')$ as 
follows:
 Given automorphisms $\alpha \in \mbox{ParAut}(E'')$, $\beta \in
\mbox{ParAut}(E')$ and a parabolic extension $[E]=(E,i,j) \in \mbox{ParExt}(E'',E')$  
we define the parabolic extension $\beta[E]\alpha$ to be the extension 
$(E,\beta \circ i,j \circ \alpha)$.
 \vs     
Now fix a \pa ~\ex ~
$[E]$ of $E''$ by $E'$.
The proof of the following lemma is analogous to the corresponding
statement for ordinary vector bundles.
\vsp
{\bf Lemma \, 2.32.}\, {\it
{\bf (a)} The orbit of $[E]$ under this action is the set of
equivalence class of \pa~  \ex s which have their middle terms
isomorphic to $E$ as \pa ~bundles.

{\bf (b)} The stablizer of $[E]$ under this action is precisely
the subgroup of  \\${\rm ParAut}(E'')\times {\rm ParAut(E')}$
consisting of elements of the form $(\alpha,\beta)$  such that there 
exists a \pa ~automorphism of $E$ which takes $E'$ to itself and 
induces $\alpha$ on $E''$ and $\beta$ on $E'$.
}
\vsp
For the convenience of the reader we summerize below, in one place the
 notations used in this paper. Some of the notations are introduced
above while the rest will be introduced subsequently.

\medskip
\centerline{\bf Summary of notation}
\medskip
\begin{supertabular}{ll}
%
%
%
$X$         & $=$ a smooth projective geometrically irreducible curve\\
            & ~~~over the finite field $\F_q$. \\

$\bar{X}$   & $=$ the curve  $X \otimes_{ \sF_q}{\bar{\F}_q}$,\\
            & ~~~where ${\bar{\F}_q}$ is an algebraic closure of $\F_q$.\\

$Z_X(t)$    & $=$ the zeta function of the curve $X$.\\

$X_{\nu}$   & $= X \otimes_{ \sF_q }{ \F_{q^{ \nu}}}$, where 
                $\F_{q^{ \nu}} \subset {\bar{\F}_q}$ is a finite\\
            & ~~~field with $q^{\nu}$ elements.\\[2mm]
 
            & ~~~For positive integers $n$ and $m$ and non-negative \\ 
            & ~~~integers $r_1,\ldots ,r_m$ with $r_1+\ldots +r_m =n$, \\
$Flag(n,m,r_j)$
            & $=$ the variety of all flags $k^n=F_1\supset\ldots\supset$\\
            & ~~~$F_m\supset F_{m+1}=0$ of vector subspaces in $k^n$,\\
            & ~~~with $dim(F_j/F_{j+1}) = r_j$.\\

$|J(\F_q)|$ 
            & $=$ the number of $\F_q$-rational points of the \\
            & ~~~Jacobian of $X$.\\

$S$         & $ =$ a finite set of $k$-rational points of $X$.\\
            & ~~~(These are the parabolic vertices.)\\

$m_P$       & $=$ a fixed positive integer defined for each $P \in S$.\\[2mm]

            & ~~~For $P\in S$, and $1 \le i  \le m_P$,\\
$\alpha$    & $ =( \alpha^P_i )$ is the set of allowed weights. \\ [2mm] 

            & ~~~For $P\in S$, and $1 \le i  \le m_P$,\\
$R$         & $ =( R^P_i )$, the quasi-parabolic data (or simply 
                      `data').\\ 

$n(R)$      & $=\sum_{i=1}^{m_P}R^P_i$, the rank of the data $R$.\\

$L$         & $=$ a sub-data of $R$ and \\

$R-L$       & $=$ the complementary sub-data defined by\\ 
            & ~~~$(R-L)^P_i = R^P_i -  L^P_i $.\\   

${\cal L}$  & $=$ a line bundle on $X$.\\

$E$         & $=$ a vector bundle with a parabolic structure with\\ 
            & ~~~data $R$.\\

$J_R({\cal L})$
            & $=$ the set of isomorphism classes of quasi-parabolic\\ 
            & ~~~vector bundles with data $R$, and determinant ${\cal L}$.\\ 

$\alpha(R)$ & $= \sum_P \sum_{i=1}^{m_P}{R^P_i \alpha^P_i}$, the
              parabolic contribution \\
            & ~~~to the degree.\\

$\mbox{deg}(E)$
            & $=$ the ordinary degree of $E$.\\            

$\mbox{pardeg}(E)$
            & $= \mbox{deg} (E) + \alpha (R)$, the parabolic degree of $E$.\\ 

$\mbox{par}\mu(E)$ 
            & $ = \mbox{pardeg} (E)/\mbox{rank}(E)$, the parabolic\\
            & ~~~slope of $E$.\\[2mm]
            
            & ~~~For $P\in S$, $1 \le i  \le m_P$, and $1 \le k \le r$,\\
$I$         & $=(I^P_{i,k})$, the intersection type of Nitsure, which is\\
            & ~~~a partition of $R$.\\ 
$\ell (I)$  & $=$the length of the intersection type $I$. \\[2mm]
 
            & ~~~For $j \le r$, we have\\
$R^I_j$     & $=$ the sub-data defined by $ (R^I_{j})^P_i = I^P_{i,j}$,\\

$R^I_{\le j}$
            & $=$ the sub-data defined by $(R^I_{\le j})^P_i = 
                    \sum_{k \le j}I^P_{i,k},$\\

$R^I_{\ge j}$
            & $=$ the sub-data defined by $(R^I_{\ge j})^P_i= 
                   \sum_{k \ge j}{I^P_{i,k}}.$\\

$I_{\le j}$ & $=$ the partition of $R^I_{\le j}$ defined by 
                $(I_{\le j})^P_{i,k}=I^P_{i,k}$\\
            & ~~~where $k \le j$.\\
$I_{\ge j}$ & $=$ the partition of $R^I_{\ge j}$ defined by $I$
                $(I_{\ge j})^P_{i,k}=I^P_{i,k}$\\
            & ~~~where $k \ge j$.\\   
${\cal M}_{R,{\cal L}}$ 
            & $=$ the moduli space of parabolic semistable bundles\\  
            & ~~~with the data $R$ and determinant ${\cal L}$.\\

${\cal M}_{R,{\cal L}}^s$ 
            & $=$ the open sub variety of ${\cal M}_{R,{\cal L}}$\\
            & ~~~corresponding to the parabolic stable bundles.\\[2mm]

            & ~~~For parabolic bundles $E', E$ and $E''$ with\\ 
            & ~~~data $R', R$ and $R''$, we denote by\\
$[E]$       & $ = (E,i,j)$, a parabolic extension of $E''$ by $E'$.\\

$\mbox{ParExt}(E'',E')$
            & ~~~the set of equivalence classes of parabolic \\
            & ~~~extensions of $E''$ by $E'$.\\

$\beta_{R}({\cal L})$ 
            & $= \sum(1/|\mbox{ParAut}(E)|)$\\
            & ~~~where summation is over all $E \in J_R({\cal L})$ such \\
            & ~~~that $E$ is parabolic semistable.\\    

$J_R({\cal L},I)$
            & $=$ the set of isomorphism classes of parabolic \\
            & ~~~bundles with weights $\alpha$, of intersection\\
            & ~~~type $I$, and determinant ${\cal L}$.\\

$\beta_R({\cal L},I)$
            & $ = \sum (1/|\mbox{ParAut}(E)|)$,\\
            & ~~~where the summation is over all $E$ in \\
            & ~~~$J_R({\cal L},I)$.\\[2mm]
            
            & ~~~We also write \\ 
$\beta_R(d,I)$
            & $=\beta_R({\cal L},I)$, since $\beta_R({\cal L},I)$ \\
            & ~~~depends on ${\cal L}$ only via its degree $d=deg({\cal L})$.\\

${\cal F}_R$
            & $\displaystyle= \prod_{P \in S} \mbox{Flag}(n(R),m_P,R^P_i)$\\

$f_R(q)$    & $=$ the number of $\F_q$-valued points of the variety 
               ${\cal F}_R$.\\   

$C(I;~ d_{1},\ldots, d_{r})$
            & $=$ the integer defined by equation (3.8).\\ 

$\sigma_k(I)$
            & $\displaystyle=\sum_{P\in S}\sum_{i>t}\sum_{l<r-k+1}
              {I^P_{i,r-k+1}I^P_{t,l}}$.\\

$\sigma_R(I)$
            & $\displaystyle= \sum_k{\sigma_k(I)}$.\\[2mm]
\end{supertabular}

\pagebreak

\begin{supertabular}{ll}
            
            & ~~~For a vector bundle $F$\\

$\chi(F)$   & $=$ the Euler characteristic.\\ 
$\chi{\scriptstyle
\left( 
\begin{array}{ccc}
{\nu_1} & {\dots} & {\nu_r} \\ {\delta_1} & {\dots} & {\delta_r}
\end{array}
\right)}$   & $=$ the numerical function of Desale-Ramanan \\
            & ~~~defined by the equation (3.20).\\

$\underset{\circ }{ \sum }$ 
            & ~~~denotes the summation over all \\
            & ~~~$(d_{1},\ldots,d_{r})\in \Z^r$ with $\sum_i{d_i}=d$\\
            & ~~~and satisfying equation (3.7).\\ 
            & \\
$\tau_{n(R)}(q)$
            & $\displaystyle =\frac{q^{(n(R)^2-1)(g-1)}}{q-1}Z_X(q^{-2})\ldots 
              Z_X(q^{-n(R)})$.\\

$\tilde {f}_R(t)$
            & $=$ the rational function corresponding to $f_R$\\ 
            & ~~~given by the equation (3.30).\\

$\tilde {\tau}_{n(R)}(t)$
            & $=$ the rational function corresponding to $\tau_{n(R)}$\\ 
            & ~~~given by the equation (3.31).\\
           
$Q_{R,d}(t)$
            & $=t^{n(R)^2(g-1)}(1+t^{-1})^{2g}\tilde{\beta}_R(d)$,\\ 
            & ~~~this is the main function for the recursion.\\

$Q_R(t)$    & $\displaystyle=t^{n(R)^2(g-1)}\tilde {f}_R(t)\tilde
               {\tau}_{n(R)}(t)$.\\  

$P_{R,d}$   & $=$ the power series whose coefficients compute the\\
            & ~~~Betti numbers of the moduli space of parabolic \\
            & ~~~stable bundles with data $R$ and degree $d$.\\ 

$N_R(I;~ d_{1},\ldots, d_{r})$
            & $=$ the integer given by the formula (3.38).\\
                                                         
$Y$         & $=$ a smooth projective variety over $\F_q$. \\


$N_{\nu}$   & $=$ the number of $\F_{q^{\nu}}$-rational points of $Y$.\\[2mm]
           
            &~~~For $i= 1,\ldots,2g$, we have \\  
$\omega_i$  & $=$ a fixed algebraic integer of norm $q^{1/2}$. \\  
$h(u,v_{1},\ldots,v_{2g})$ 
            & $=$ a rational function given by the equation (4.2)\\

$p(u,v_{1},\ldots,v_{2g})$
            & $=$  the numerator occuring in the equation (4.2).\\

$(a_{J.j})$
            & $=$ the coefficients occuring in (4.3).\\

$J$         & $=$ the multi-index $J=(i_{1},i_{2},\ldots, i_{2g})$,\\ 

$|J|$       & $= \sum_{r=1}^{2g}{i_{r}}$, and \\

$v^{J}$     & $= v_{1}^{i_{1}}v_{2}^{i_{2}}\ldots  v_{2g}^{i_{2g}}$.\\

$N$         & $=$ the `weighted degree' of $p(u,v_{1},v_{2},\ldots,v_{2g})$.\\

$b_{J,j}$   & $=$ the coeffiients of $h$ defined in (4.8).\\

$f_{\ge 0}(u,v_{1},\ldots,v_{2g})$ 
            & $=$ the function defined by the equation (4.13).\\ 

$M_r$       & $=f_{\ge0}(q^{r},\omega_{1}^{r},\ldots,\omega_{2g}^{r})$.\\

$Z_1(t)$    & $=$ the formal power series defined in (4.14).\\

$Z_2(t)$    & $=$ the formal power series defined in (4.15).\\

$Z(t)$      & $= Z_1(t)Z_2(t)$.\\[2mm]
            
            &~~~For a meromorphic function $h$ on a disc \\
            &~~~in $\C$, and $\alpha >0$,\\
$\mu(h,\alpha)$
            & $=$ the number of zeros minus the number of poles \\
            &~~~counted with multiplicities of $h$ with norm $\alpha$.\\

$P(T)$      & $=$ the polynomial defined by the equation (4.19).\\ 

$ M'_R(I;~d)$
            & $=$ the integer given by the formula (5.4).\\[2mm]
            
            & ~~~For a real number $\lambda$,\\
$M_R(I;~\lambda)$
            & $=$ the integer given by the formula (5.5).\\

$Q_{R,d}^{\lambda}(t)$
            & $=$ the rational function defined in (5.7).\\

$S_{R,d}^{\lambda}(t)$
            & $=$ the rational function defined in (5.8).\\

$\underset{{\circ_{\lambda}}}{ \sum }$ 
            & ~~~denotes the summation over 
              $(d_{1},\ldots,d_{r}) \in \Z^r$\\
            & ~~~such that  $\sum_i{d_i}=d$ and the equation (5.9) holds.\\

$Q_{R,d}^{\lambda^-}$ 
            & $=Q_{R,d}^{\lambda-\epsilon}$ for $\epsilon$ small enough \\
            & ~~~such that the function $Q_{R,d}^{\lambda}$ has no jumps \\
            & ~~~in the interval $[\lambda-\epsilon, \lambda)$.\\ 

$S_{R,d}^{\lambda^-}$ 
            & $=S_{R,d}^{\lambda-\epsilon}$ for $\epsilon$ small enough \\
            & ~~~such that the function $S_{R,d}^{\lambda}$ has no jumps \\
            & ~~~in the interval $[\lambda-\epsilon, \lambda)$.\\ 

$\Delta Q_{R,d}^{\lambda}$
            & $= Q_{R,d}^{\lambda}-Q_{R,d}^{\lambda^-}$.\\

$\Delta S_{R,d}^{\lambda}$
            & $=S_{R,d}^{\lambda}-S_{R,d}^{\lambda^-}$.\\

$\delta_R(L)$ 
            & $=$ the integer given by the equation (5.10).\\

$d(\la,L)$
            & $\displaystyle= n(L) \la - \alpha(L)$. \\

$g_R(I;~d)$
            & $=$ the rational function given by (5.17).\\

$\sigma'_R(I)$
            & $\displaystyle= \sum_{P\in S}\sum_{k>l,i<t}
               {I^P_{i,k}I^P_{t,l}}$.\\

$M_g(I;~\lambda)$ 
            & $=$ the integer given by the formula (5.25).\\

$P_R(t)$    & $=$ the rational function(polynomial) defined by \\
            & ~~~the equation (5.26).\\[2mm] 

            & ~~~For a data $R$ with rank $n(R)=2$, \\
$T$         & $=$ the subset of $S$ consisting of parabolic vertices 
              where\\       
            & ~~~the parabolic filtration is non-trivial.\\

$T_I$       & $=~~ \{P \in T|I^P_{1,1} =0\}$.\\

$\chi_I$    & $=$ a characteristic function on $T$, defined by \\
            & ~~~the equation (6.3).\\  

$\psi_I$    & $= \sum_{P \in T}\chi_I(P)(\alpha^P_1 -\alpha^P_2)$.\\

$a_I$       & $= 1$ if $d+ [\psi_I]$ is even, and\\

            & $= 0$ if $d+ [\psi_I]$ is odd.\\

$\delta^P$  & $ =\alpha^P_1 -\alpha^P_2$.\\
\end{supertabular}  
\medskip
\section {The inductive formula}
\medskip
 \centerline{\bf The use of Parabolic-Harder-Narasimhan Intersection types}
   \medskip

 In this section we use the quasi-parabolic Siegel
 formula to obtain a recursive formula for the Poincar\'e polynomial of
 the moduli space of \psb s  when the condition `par semi-stable = par stable' holds.
 \vs    
       The left hand side of the quasi-parabolic Siegel formula (2.19) can
       be split into the summations coming from the parabolic semistable
bundles and the unstable ones. In view of this we first define
 $$\beta_{R}({\cal L}) = \sum{\frac {1}{|\mbox{ParAut}(E)|}}
   \eqno(3.1)$$ 
where summation is over all $E \in J_R({\cal L})$ such that $E$ is 
parabolic semistable.
\vs
  We assume that the data $R$ and degree $d$ are so chosen
that the condition `par semi-stable = par stable' holds.
In particular by lemma (2.22) this implies that for any such \pssb
~$E$, $|ParAut(E)| = q-1$. Hence
$$|{\cal M}_{R,{\cal L}}(\F_q)| = (q-1)\beta_R({\cal L}) \eqno(3.2)$$
is the number of $\F_q$- rational points of the moduli space
of parabolic \\ semistable bundles with the data $R$ and determinant 
${\cal L}$.
 \vs
   Now we have to take care of the unstable part of the
   summation (2.19). This summation can be further split into \phn
   ~intersection types. For these considerations, we make the
   following definitions:
  \vs     
Let $I$ be a partition of $R$ with $\ell (I)=r$. Let $J_R({\cal
L},I)$ denote the set of isomorphism classes of parabolic bundles
with data $R$, of intersection type $I$, and determinant ${\cal L}$.
\vs
Let $$\beta_R({\cal L},I) = \sum {\frac{1}{|\mbox{ParAut}(E)|}} \eqno(3.3) $$
where the summation is over all $E$ in $J_R({\cal L},I)$. Note that   
$\beta_R({\cal L},I)=\beta_{R}({\cal L})$ for the unique $I$ which has $\ell(I)=1$.
 \vs      
The summations occurring in (3.1) and (3.3) are finite
because the \pa~ \\bundles  of fixed intersection type form a bounded
family, so it is dominated by a variety, hence has only finitely many 
$\F_q$-rational points.
 \vs       
      
Now the quasi-parabolic Siegel formula (2.19) can be restated as
$$\sum_{r \ge 1}\sum_{\{I|\ell (I)=r\}}{\beta_R({\cal L},I)}
   =  \frac {f_R(q)q^{(n(R)^2-1)(g-1)}}{q-1}Z_X(q^{-2})\ldots Z_X(q^{-n(R)}) 
\eqno(3.4)$$
  where $f_R(q)$ is given by (2.17).
\vs

\medskip
\centerline{\bf Computation of the function $\beta_R({\cal L},I)$}
\medskip
The main step in the induction formula is to use the parabolic Harder 
-Narasimhan filtration to give a formula for $\beta_R({\cal L},I)$ when
 $\ell(I)>1$, in terms of
 $\beta_{R'}({\cal L}')$ of lower rank bundles. This we do in the  
following proposition which is an analogue of proposition (1.7) of Desale
and Ramanan[D-R].
\vsp
{\bf Proposition \, 3.5.}\, {\it 
{\bf (a)} The numbers $\beta_R({\cal L},I)$ and $\beta_R({\cal L})$
depend on ${\cal L}$ only via its degree $d=deg({\cal L})$
(hence they can be written as $\beta_R(d,I)$ and $\beta_R(d)$ resp.).

{\bf (b)} $\beta_R(d,I)$ satisfies the following recursive relation
$$ \beta_R(d,I) = \sum_{\circ } q^{C(I;~d_{1},\ldots, 
 d_{r})} |J(\F_q)|^{r-1}
 \prod_{k=1}^r{\beta_{R^I_k}(d_k)}              \eqno(3.6)$$

where $\underset{\circ}{\sum }$ denotes the summation over all 
$(d_{1},\ldots,d_{r}) \in \Z^r$
with $\sum_i{d_i}=d$ and  satisfying the following inequalities
$$
\frac {d_1 + \alpha(R^I_1)}{n(R^I_1)} > \frac {d_2 + \alpha(R^I_2)}{n(R^I_2)}
> \ldots >\frac {d_r + \alpha(R^I_r)}{n(R^I_r)} \eqno(3.7)
$$
Here $|J(\F_q)|$ denotes the number of
    $\F_q$-valued points of the Jacobian of $X$, and
$$
\begin{array}{ccl}
C(I;~ d_{1},\ldots, d_{r})& = \sum_{P\in S} 
    \sum_{k>l,i>t}{I^P_{i,k}I^P_{t,l}} & -
    \sum_{k>l}{(d_{l}n(R^I_k)-d_{k}n(R^I_k)
)}\\ & & \\ & & + \sum_{k>l}{n(R^I_l)n(R^I_k)(g-1)}
\end{array} \eqno(3.8)$$
}
\vsp
{\noindent {\bf Proof.}} We prove both parts ((a) and (b)) of the
proposition simultaneously by
induction on $n=n(R)$. If $\ell (I)=1$ then there is nothing to
prove. 
\vs
Consider a parabolic bundle $E$ with data $R$,
admitting the\\ parabolic Harder-Narasimhan filtration  
$ 0 \subset G_1 \subset\ldots \subset G_r = E$ of length $r \ge 2$.
 Let $M$ be the quotient $E/G_1$. If we give the induced parabolic
structure to $M$ then it has the data $R^I_{\ge 2}$. 
 \vs  
  Let $T$ be the set of equivalence classes of parabolic
extensions of $M$ by $G_1$ which has the property that the middle term
is isomorphic to $E$ as a parabolic bundle. 
     By lemma (2.32 (a)) $T$ is same as the orbit of $[E]$ under the action of 
 $\mbox{ParAut}(M) \times  \mbox{ParAut}(G_1)$ on $
\mbox{ParExt}(M,G_1)$, 
hence
         $$|T| = \frac{|\mbox{ParAut}(M)|
             |\mbox{ParAut}(G_1)|}{|\mbox{stabilizer of}~
             [E]|}. \eqno(3.9)$$
 
Note that every parabolic automorphism of $E$ takes $G_1$ to
             itself (hence also $M$).
This implies that we get a group homomorphism 
$$\mbox{ParAut}(E) \stackrel{\phi}{\longrightarrow} 
       \mbox{ParAut}(G_1) \times  \mbox{ParAut}(M). \eqno(3.10)$$
Now by lemma (2.32(b)) the stabilizer of $[E]$ is the
   image of $\phi$, while the kernel of $\phi$ is equal to 
$I + H^0(X,{\mbox{${\cal P}ar{\cal 
             H}om~$}}(M,G_1))$.

Combining all this  we get
$$|T| = \frac {|\mbox{ParAut}(M)| |\mbox{ParAut}(G_1)|
|\mbox{ParHom}(M,G_1)|}{|\mbox{ParAut}(E)|}.  \eqno(3.11)$$

     By definition 
      $$ \beta_R({\cal L},I) =\sum_{E\in
J_R({\cal L},I)}{\frac {1}{|\mbox{ParAut}(E)|}}
   \eqno(3.12)$$  
          which is $$\sum_{(M,G_1)}\sum_{{\cal E}} {\frac
{1}{|\mbox{ParAut}(E)||T|}} \eqno(3.13)$$
where the first summation extends over all pairs $(M,G_1)$ with 
$G_1$, a \pssb ~with data
$R^I_1$, and $M$, parabolic bundle with data $R^I_{\ge 2}$ and 
 intersection type $I_{\ge 2}$,  
 such that $\mbox{det}(M)\otimes \mbox{det}(G_1)= \mbox{det}(E)$. The 
second summation extends over the set ${\cal E} = \mbox{ParExt}(M,G_1)$.
   By (3.11), the right hand side of the above expression (3.13) reduces to
 $$ 
\sum_{(M,G_1)}{\frac{1}{|\mbox{ParAut}(M)| |\mbox{ParAut}(G_1)| 
q^{\chi({\mbox{${\cal P}ar{\cal H}om\,$}}(M,G_1))}}},   \eqno(3.14)$$
where $$\chi({\mbox{${\cal P}ar{\cal H}om\,$}}(M,G_1)) =
\mbox{dim}_{\sF_q}(\mbox{ParHom}(M,G_1))-\mbox{dim}_{\sF_q}(\mbox{ParExt}(M,G_1))
\eqno(3.15)$$
is the Euler characteristic of the sheaf ${\mbox{${\cal P}ar{\cal H}om\,$}}(M,G_1)$.
\vs
We define certain numerical functions which depends only on the
partition $I$ as follows:
 $$\sigma_k(I)= \sum_{P\in S}\sum_{i>t}\sum_{l<r-k+1}
{I^P_{i,r-k+1}I^P_{t,l}} ~~~\mbox{and}~~~ \sigma_R(I) = \sum_k{\sigma_k(I)} \eqno(3.16)$$
   
     Then it can be checked that $\sigma_1(I)$ is the length of the torsion
 sheaf  ${\cal S}_1(I)$,
 which is defined by the following exact sequence:
 $$0 \longrightarrow {\mbox{${\cal P}ar{\cal H}om\,$}}(M,G_1)
\longrightarrow \shom(M,G_1) \longrightarrow {\cal S}_1(I) \longrightarrow
0.          \eqno(3.17) $$

Using the fact that
$$\chi({\mbox{${\cal
             P}ar{\cal H}om\,$}}(M,G_1)) =  \chi(M^* \otimes G_1) -
             \sigma_1(I),
   \eqno(3.18)$$ 

the sum (3.14) becomes 
          $$ = \sum_{M,G_1}{\frac{q^{\sigma_1(I)}}{|\mbox{ParAut}(M)| 
|\mbox{ParAut}(G_1)|q^{\chi(M^*\otimes G_1)}}}.       \eqno(3.19) $$

Recall that Desale-Ramanan[D-R] introduced certain numerical functions 
 $$\chi
{\scriptstyle
\left( 
\begin{array}{ccc}
{\nu_1} & {\dots} & {\nu_r} \\ {\delta_1} & {\dots} & {\delta_r}
\end{array} 
\right)} = {{\sum_{k>l}({\delta_{l}\nu_k-\delta_{k}\nu_k})}} +
{\sum_{k>l}{\nu_l\nu_k(g-1)}}.     \eqno(3.20) $$ 
 With this definition of $\chi$, we have the following equality 

$$\chi
{\scriptstyle 
\left( 
\begin{array}{cc}
{n(R^I_1)} & {n(R)-n(R^I_1)} \\ {d_1} & {d-d_1}
\end{array} 
\right)} = \chi(M^* \otimes
G_1)         \eqno(3.21)$$
as in [D-R].
Now by (3.21), the  sum (3.19) equals 
$$ \sum_{d_1}{{q^{\sigma_1(I)-\chi
{\scriptstyle
\left( 
\begin{array}{cc}
{n(R^I_1)} & {n(R)-n(R^I_1)} \\ {d_1} & {d-d_1}
 \end{array} 
\right)}}}\sum_{(\eta,\gamma)}
{\sum_M \frac{1}{\mbox{ParAut}(M)}\sum_{G_1}{\frac{1}{\mbox{ParAut}
(G_1)}}}},
       \eqno(3.22) $$ 
where the first summation in (3.22) is over all integers $d_1$ with 
     $$  (d_1 +\alpha(R^I_1))/n(R^I_1) >
(d-d_1 + \alpha(R^I_{\ge 2}))/(n-n(R^I_1)).   \eqno(3.23) $$

The second summation in (3.22) is over isomorphism classes of
 line bundles $\eta$ and $\gamma$ such that
$\eta \otimes \gamma = {\cal L}$. The
third one is over all parabolic bundles $M$ with data $R$,
 having intersection type $I_{\ge 2}$, and determinant $\eta$. The fourth 
summation in (3.22) is
     over all semi-stable parabolic bundles $G_1$ with data $R^I_1$,
    and determinant $\gamma$.
  This expression is equal to $$ \sum q^{\sigma_1(I)-\chi
{\scriptstyle
\left(\begin{array}{cc}{n(R^I_1)} & {n-n(R^I_1)}\\ {d_1} &
{d-d_1}\end{array}\right)}}\sum
{\beta_{R^I_{\ge 2}}(\eta,I_{\ge 2})\beta_{R^I_1}(\gamma)}.  
\eqno(3.24)$$
 Now note that by induction, the terms inside the summation are
     independent of ${\cal L}$, hence part (a) of the proposition
 follows. 
 From now on we write $\beta_R(d)$ and $\beta_R(d,I)$ for  
$\beta_R({\cal L})$ and $\beta_R({\cal L},I)$. 
\vs
By Desale-Ramanan [D-R] we have  the relation 
 $${\scriptstyle \chi
\left( 
\begin{array}{cc}
{n(R^I_1)} & {n-n(R^I_1)} \\ {d_1} & {d-d_1}
\end{array} 
\right) 
  + \chi
\left( 
\begin{array}{ccc}
{n(R^I_2)} & {\dots} & {n(R^I_r)} \\ {d_2} & {\dots} & {d_r}
\end{array} 
\right) 
 = \chi
\left( 
\begin{array}{ccc}
{n(R^I_1)} & {\dots} & {n(R^I_r)} \\ {d_1} & {\dots} & {d_r}
\end{array} 
\right)}     \eqno(3.25)$$ 
 Using this and the induction hypothesis for
$\beta_{R^I_{\ge 2}}(d-d_1,I_{\ge 2})$ we obtain the following equality:
$$\beta_R(d,I) = \sum q^{\sigma_R(I)- \chi
{\scriptstyle 
\left( 
\begin{array}{ccc}
{n(R^I_1)} & {\dots} & {n(R^I_r)} \\ {d_1} & {\dots} & {d_r}
\end{array} 
\right)}} |J_{\sF_q}|^{r-1}
\prod_{k=1}^r{\beta_{R^I_k}(d_k)}
   \eqno(3.26) $$
 As 
$$C(I;~ d_{1},\ldots, 
 d_{r}) = \sigma_R(I) - \chi
{\scriptstyle 
\left( 
\begin{array}{ccc}
{n(R^I_1)} & {\dots} & {n(R^I_r)} \\ {d_1} & {\dots} & {d_r}
\end{array} 
\right)}       \eqno(3.27) $$
the proof of the proposition is complete. $\hfill{\Box}$

\medskip
\centerline{\bf The recursive formula}
\medskip
The inductive expression for  $\beta_R(d)$ can now be written
as
 $$ f_R(q)\tau_{n(R)}(q) - \sum_{r\ge 2}\sum_{\{I|\ell (I)=r\}}
\sum_{\circ} q^{C(I;~ d_{1},\ldots, d_{r})} |J_{\sF_q}|^{r-1} 
   \prod_{k=1}^r{\beta_{R^I_k}(d_k)}
    \eqno(3.28)$$
where $$\tau_{n(R)}(q) =
 \frac{q^{(n(R)^2-1)(g-1)}}{q-1}Z_X(q^{-2})\ldots 
Z_X(q^{-n(R)}).       \eqno(3.29)$$

    We now base change from  $\F_q$ to  $\F_{q^{\nu}}$.
For the curve $X_{\nu}$ defined in the section 2,
the  $\beta_R(d,q^{\nu})$ will be a function of 
$q^{\nu}$ and $\omega_i^{\nu}$ for $i= 1 ,\ldots,2g$.
 \vs 
In the light of the induction formula
and equation (2.17)
  we get that the function $\beta_R(d,q^{\nu})$ is a polynomial in 
    $ \omega_i^{\nu}$ for $i = 1 ,\ldots, 2g$, and is a rational
function in $q^{\nu}$ with the property that the denominator has
factors only of the form $q^{\nu n_0} (q^{\nu n_1}-1)(q^{\nu n_2}-1)\ldots
(q^{\nu n_k}-1)$, with $n_i\ \ge 1 ~\mbox{for}~ i \ge 1$.
   For such a functions, one can substitute
 $-t^{-1}$ for  $\omega_i$ and $t^{-2}$ for $q$, to obtain a 
 new rational function. We denote this operation by $\phi
 \rightarrow  \tilde {\phi}$.
 For example $$\tilde {f}_R(t) =  \frac {t^{-2 \dim{\cal F}_R}
\prod_{i=1}^{n(R)}{(1-t^{2i})}^{|S|}}
{\prod_{P\in S}\prod_{\{i| {R^P_i} \ne 0\}}\prod_{l=1}^{R^P_i}{(1-t^{2l})}}
\eqno(3.30)$$ 
 and   $\tilde {\tau}_{n(R)}(t)$  can be computed to be
$$\frac{t^{-2n(R)^2(g-1)}\prod_{i=1}^{n(R)}{(1+t^{2i-1})^{2g}}}
{(1-t^{2n(R)})\prod_{i=1}^{n(R)-1}{(1-t^{2i})^2}}
    \eqno(3.31) $$ 
This substitution is an important step in the computation of the
Poincar\'e polynomial for the moduli space  because of the proposition
(4.34) of the next section.

Now we shall define rational functions $Q_{R,d}(t)$ and $Q_R(t)$ by
$$Q_{R,d}(t) =
 t^{n(R)^2(g-1)}(1+t^{-1})^{2g}\tilde{\beta}_R(d)    \eqno(3.32)$$ 
and 
$$Q_R(t)=t^{n(R)^2(g-1)}\tilde {f}_R(t)\tilde {\tau}_{n(R)}(t).   
  \eqno(3.33)$$

Observe that if the  condition `par semi-stable = par stable' is satisfied and
 if we define 
$$P_{R,d}(t) = t^{ 2 \dim{\cal F}_R+2(n(R)^2-1)(g-1)}(t^{-2}-1)\tilde{\beta}_R(d) 
 \eqno(3.34)$$
 then we have the relation $$P_{R,d}(t) = \frac{
 t^{2 \dim{\cal F}+n(R)^2(g-1)}(1-t^2)}{(1+t)^{2g}}Q_{R,d}(t).  \eqno(3.35)$$  

 Now by proposition (4.34) and the fact that the dimension of the
 moduli space of \pssb s is equal to $\dim{\cal F}_R+(n(R)^2-1)(g-1)$,
 we get that $P_{R,d}$ is a power
 series in $t$ which computes the Betti numbers of the 
 moduli space of \psb s with the given data $R$ and degree $d$.
 \vs
    If we perform the tilde operation on the original formula, we get the
 following recursive formula.
\vsp
{\bf Theorem \, 3.36.}\, {\it  
The functions $Q_{R,d}$ and $Q_R$ defined by {\rm (3.32)} and
 {\rm (3.33)} satisfy the following recursion formula.
 $$Q_R(t) =  \sum_{r \ge 1}\sum_{\{I|\ell (I)=r\}} \sum_{\circ} t^{2N_R(I;~ d_{1},\ldots, 
 d_{r})} \prod_{k=1}^r{Q_{R^I_k,d_k}(t)}     \eqno(3.37)$$
   where the second summation extends over all partitions $I$ of $R$
 of length $r$, and where 
     $$N_R(I;~ d_{1},\ldots, d_{r}) =
    \sum_{k>l}{(d_{l}n(R^I_k)-d_{k}n(R^I_l))}-\sum_{P\in S} 
    \sum_{k>l,i>t}{I^P_{i,k}I^P_{t,l}} . \eqno(3.38)$$
}

\section{The substitution $\omega_i \rightarrow -t^{-1},~q \rightarrow
t^{-2}$}
In this section, we justify the substitution $\omega_i \rightarrow
-t^{-1} ~~\mbox{and}~~~q \rightarrow t^{-2}$, which gives us a recipe to
compute the Poincar\'e polynomial of the moduli spaces, directly from
the computation of the $\F_q$-rational points.
  \vs
  This substitution was briefly sketched in [H-N] for the
rational function which counted the $\F_q$ rational points of the
moduli space of stable bundles when rank and degree are coprime.
We formulate and prove this in a more general setup, which we have used
in the body of the paper.
 \vs
Let $Y$ be a smooth projective variety over $\F_q$. Let $N_{\nu} =
|Y(\F_{q^{\nu}})|$ and let $\omega_1,\ldots ,\omega_{2g}$ be fixed algebraic
integers of norm $q^{1/2}$. 
   Our basic assumption is that $N_{\nu}$ is given by some formula
     $$N_{\nu} = h(q^{r},\omega_{1}^{r},\ldots,\omega_{2g}^r)\eqno(4.1)$$
where $h(u,v_{1},\ldots,v_{2g})$ is a rational function  
of the form
$$\frac {p(u,v_{1},\ldots,v_{2g})}{u^{n_0}(u^{n_1}-1)\ldots
(u^{n_k}-1)}\eqno(4.2)$$ 
where  $p(u,v_{1},\ldots,v_{2g}) \\ \in \Z[u,v_{1},\ldots,v_{2g}]$
is a polynomial with integral coefficients, and where 
 $n_i \ge 1$ for all $i > 0$ and $n_0 \ge 0$. 
We wish to write down the Poincar\'e polynomial of
$Y$ in terms of the function $h$.

We first write down the function $h$ as a suitable series and bound
the coefficients.
We can expand the numerator occurring in the expression for $h$ as 
$$p(u,v_{1},\ldots,v_{2g})= \sum_{l=0}^{N} 
\sum_{|J|+2j=l}{a_{J,j}v^{J}u^{j}} \eqno(4.3)$$ 
where $J$  
denotes the multi-index $J=(i_{1},\ldots, i_{2g})$, 
$|J| = \sum_{r=1}^{2g}{i_{r}}$, 
and $v^{J} = v_{1}^{i_{1}}\ldots  v_{2g}^{i_{2g}}$.
Let $C>0$ be any fixed integer such that 
$|a_{J,j}| < C~\mbox{\rm{for all $J$, $j$}}.$
The integer $N$ in the summation above can be taken to be
the `weighted degree' of $p(u,v_{1},\ldots,v_{2g})$ where the variable
$u$ is given weight $2$.

\vs
We can rewrite $h$ as 
$$\frac{1}{u^n(1-u^{-n_1})\ldots (1-u^{-n_k})}\sum_{l=0}^{N}
\sum_{|J|+2j=l}{a_{J,j}v^{J}u^{j}}  \eqno(4.4)$$ 
where $n= \sum_{i=0}^{k}{n_i}$. 
Expanding each $1/(1-u^{-n_i})$ as a power series in $u_{-1}$, we get
$$h = \frac{1}{u^n}\sum_{l=0}^{N}
\sum_{|J|+2j'=l}\sum_{i\le 0}{a_{J,j'}b_iv^{J}u^{j'+i}} \eqno(4.5)$$
where $b_i$ is the cardinality of the set of k-tuples of non-negative
integers $(a_1,\ldots, a_k)$ such that $\sum_{r=1}^{k}{a_rn(R)} = -i.$
Clearly, we have
$$b_i \le (-i+1)^k  \eqno(4.6)$$
Now the right hand side of the equation (4.5) becomes    
$$ \frac{1}{u^n}\sum_{l \le N} 
\sum_{|J|+2j=l}{ \left( \sum_{j'+i=j}{{a_{J,j}b_i}}\right) v^{J}u^{j}}. \eqno(4.7)$$
Define 
$$b_{J,j}= \sum_{j'+i=j} a_{J,j}b_i   \eqno(4.8)$$
This is a finite sum, which 
makes sense for every $J,j$ such that $|J|+2j \le N$. 
In terms of these $b_{J,j}$, 
the expression for $h$ can be written as 
$$ h = \frac{1}{u^n}\sum_{l\le N}
\sum_{|J|+2j=l}{b_{J,j}v^{J}u^{j}}. \eqno(4.9)$$
Note that in the above series, there are only finitely 
many positive powers of $u$ and infinitely many negative powers. 
The following lemma puts a bound on the coefficients $b_{J,j}$.
\vsp
{\bf Lemma \, 4.10.}\, {\it  
The coefficients $b_{J,j}$ as defined above satisfies the
following inequality
$$|b_{J,j}| \le {CN(N-j+1)^k}   \eqno(4.11)$$
}
\vsp
{\noindent {\bf Proof.}}One observes that
$$ |b_{J,j}| \le  \sum_{j'+i=j}|{a_{J,j}b_i}|
\le C \sum_{j'+i=j}{|b_i|}  \eqno(4.12)$$ 
where $j'$ and $i$ are as in the preceding discussion. 
In the last expression of (4.12), the number of terms is $\le N$,
and by (4.6) each term $|b_i|$ is bounded by $(-i+1)^k$. 
As $(-i+1)^k \le (N-j+1)^k$ for $j'+i=j$,
the last expression (4.12) is bounded by 
$CN(N-j+1)^k$. This proves the lemma. $\hfill{\Box}$
\vsp
 Let 
$$
\begin{array}{llll}
h_{\ge0}(u,v_{1},\ldots,v_{2g}) & = &  \sum_{l=2n}^{N}
\sum_{|J|+2j=l}{b_{J,j}v^{J}u^{j-n}} & \mbox{if}~ N \ge 2n \\
     &   &   & ~~~~~~~~~~~~~~~~~~~~\,\,(4.13)          \\
     & = & 0 & \mbox{otherwise}
\end{array}
$$ 
and $ M_r$ be
$h_{\ge0}(q^{r},\omega_{1}^{r},\omega_{2}^{r},\ldots,
\omega_{2g}^{r})$, then these numbers are well defined \\because of lemma 1.\\
Let $$Z_1(t)=\exp(\sum_{r\ge 1}{M_{r}t^{r}/r}) \eqno(4.14) $$
and $$Z_2(t)= \exp(\sum_{r\ge 1}{(N_r-M_r)t^{r}/r}),  \eqno(4.15) $$ then $
Z_1(t)$ and $Z_2(t)$ are well defined formal power series. We also define
$$Z(t)= Z_1(t)Z_2(t). \eqno(4.16) $$

       Given any meromorphic function $h$ on a disc in $ \C$, let 
$\mu(h,\alpha)$
denote the number of zeros minus the number of poles 
$h$ with norm $\alpha$, counted with multiplicities.
\vsp
{\bf Lemma \, 4.17.}\, {\it 
 {\bf (a)} $Z_2(t)$ is a non-vanishing holomorphic function on the
disc $|t|<q^{1/2}$ and $Z_1(t)$ is a rational function, hence $Z(t)$
 is a well defined meromorphic function in the region $|t|< q^{1/2}$,
 such that  $$\mu(Z(t),q^{-i/2}) = \mu(Z_1(t),q^{-i/2}).  \eqno(4.18)$$
{\bf (b)} Let $$P(T)= \sum_{i \ge 0}
(-1)^{i+1}\mu(Z(t),q^{-i/2}){T^i},   \eqno(4.19)$$ then
$$P(T)=h_{\ge0}(T^2,-T,-T,\ldots,-T). \eqno(4.20)$$
}
\vsp
{\noindent {\bf Proof.}}
     To prove that the function  $Z_2(t)$ has the above mentioned
property it is enough to verify that
the function $$g(t):= \sum_{r\ge 1}{(N_r-M_r)t^{r}/r}  \eqno(4.21) $$ is holomorphic on
the disc $|t|< q^{1/2}$. This function is $$
\sum_{l < 2n}\sum_{|J|+2j=l}{b_{J,j}\omega^{Jr}q^{jr-nr}t^{r}/r}. \eqno(4.22)$$
  The coefficient of $t^r$ is equal to
       $$\sum_{l < 2n}\sum_{|J|+2j=l}{b_{J,j}\omega^{Jr}q^{jr-nr}/r} 
 \eqno(4.23)$$     
whose modulus is bounded by  
   $$ \sum_{l < 2n}\sum_{|J|+2j=l}{|b_{J,j}|q^{r(j-n+|J|/2)}/r} \eqno(4.24) $$              
     This by lemma (4.1) is $$ \le \frac{NC}{rq^{nr}}\sum_{l < 
2n}\sum_{|J|+2j=l}{(N-j+1)^{k}q^{r(2j+|J|)/2)}}  \eqno(4.25)$$
        $$  \le  \frac{N^{2}C}{rq^{nr}} \sum_{l < 2n}
{q^{rl/2}((3N+2-l)/2)^k}   \eqno(4.26)$$   
     $$ \le \frac{N^{2}C}{2^{k}r} \sum_{l > 0}{q^{-rl/2}(3N+2-2n+l)^k}
  \eqno(4.27)$$
 which is clearly a finite sum for $r \ge 1$ because powers of $q$ decay
exponentially and the other term has polynomial growth.
Now since $(a+l) \le a^l$ for $ a \ge 2$ therefore the above summation is
bounded by
            $$ 2^{-k}N^{2}C \sum_{l > 0}{((3N+2-2n)^{k}/q^{r/2})^l}/r.
 \eqno(4.28)$$    
   Suppose $r$ is large enough such that $q^{r/2} > 2(3N+2-2n)^{k}$
the coefficient of $t^r$ has the bound
$2^{-k+1}N^{2}C(3N+2-2n)^{k}/(rq^{r/2})$
and the series
with coefficient of $t^r$ as above for large $r$ clearly has radius of
convergence $q^{r/2}$.
    Now we compute  $Z_1(t)$ as $$ \exp(\sum_{r\ge
1}\sum_{l=2n}^{N}\sum_{|J|+2j=l}{b_{J,j}\omega^{Jr}q^{jr-nr}t^{r}/r}) 
 \eqno(4.29)$$
   $$=\prod_{l=2n}^{N}\prod_{|J|+2j=l}\exp({b_{J,j}}\sum_{r\ge
1}{\omega^{Jr}q^{jr-nr}t^{r}/r})  \eqno(4.30)$$
which is equal to $$
\prod_{2n \le l \le
N}\prod_{|J|+2j=l}(1-\omega^{J}q^{j-n}t)^{(-1)b_{J,j}},
   \eqno(4.31)$$
           hence this is a rational function, and this also proves that
$Z(t)$ is a \\ meromorphic function in the region $|t| < q^{1/2}$,
and that $\mu(Z(t),q^{-i/2}) = \mu(Z_1(t),q^{-i/2})$. This
finishes the proof of part (a).
\vs
    Also from here we can read off that $$ \mu(Z_1(t),q^{-i/2}) =
(-1)\sum_{|J|+2j+2n=i}{b_{J,j}}. \eqno(4.32)$$
Clearly the polynomial $f_{\ge0}(T^2,-T,-T,\ldots,-T)$ now coincides
with $$\sum_{i\ge0}(-1)^{i+1}\mu(Z_1(t),q^{-i/2}){T^i}. \eqno(4.33)$$
Now by part (a) the proof of the lemma is complete.  $\hfill{\Box}$
\vsp
Now for the variety $Y$ if the $\F_q$-rational points are given by
 the equation (4.1) and (4.2), we get that  $Z(t)$ is the
  zeta function of $Y$ and $P(t)$ is the Poincar\'e polynomial of $Y$. 
We can restate the lemma(4.17) in terms of the Poincar\'e polynomial of
 $Y$, using Poincar\'e duality, as follows.
\vsp
{\bf Proposition \,4.34.}\, {\it 
 The function $T^{2{\rm dim}(Y)}h(T^{-2},-T^{-1},-T^{-1},\ldots,-T^{-1})$
has a formal power series expansion $\sum_{\nu \ge 0}b_{\nu}T^{\nu}$
where $b_{\nu}$ is the $\nu^{\mbox{th}}$-Betti number of $Y$ for $\nu
\le 2{\rm dim}(Y)$.
}
\section{The Closed Formula}
In this section we solve the recursion formula (theorem (3.36)) to obtain
a closed formula for the Poincar\'e polynomial of the moduli space of
parabolic stable bundles under the condition `par semi-stable = par
stable'. We do this by generalizing the method of
Zagier[Z] to the parabolic set up.
 \vs     

The induction formula can be re-written as 
 $$Q_R(x) =  \sum_{r \ge 1}\sum_{I} \sum_{\circ} 
x^{n(R)(I;~d_{1},\ldots, d_{r})}
\prod_{k=1}^r{Q_{R^I_{k},d_k}(x)}  \eqno(5.1)$$
where $x=t^2$.
The closed formula for $Q_{R,d}$ is given by the following theorem.
\vsp
{\bf Theorem \, 5.2}\, {\it 
Let $Q_{R,d}$ and $Q_R$ be formal Laurent series in
$\Q((x))$ related by the formula {\rm (5.1)}. For any $d$ and $R$
we have 
    $$ Q_{R,d}(x) = \sum_{r \ge 1}\sum_{I} 
\frac{x^{M'_R(I;~d)+M_R(I;~(d+\alpha(R))/n(R))}}{(x^{n(R^I_1)+n(R^I_2)}-1)
\ldots (x^{n(R^I_{r-1})+n(R^I_r)}-1)} \prod_{k=1}^r{Q_{R^I_{k}}(x)}  \eqno(5.3)$$
where $M'_R(I;~d)$ and $M_R(I;~\lambda)$ for a partition $I$ of $R$ and
$\lambda \in \R$ are defined by
$$ M'_R(I;~d) = -(n(R)-n(R^I_r))d - \sigma_R(I) +
  (2n(R)-n(R^I_1)-n(R^I_r))
  \eqno(5.4)$$
$$\mbox{and}~~M_R(I;~\lambda) = \sum_{k=1}^{r-1}(n(R^I_k)+n(R^I_{k+1}))
[(n(R^I_1)+\ldots + n(R^I_k))\lambda - \alpha(R^I_{\le k})].  \eqno(5.5)
$$
Here $[x]$ for a real number $x$ denotes the largest integer less than
or equal to $x$. 
}
\vsp
{\noindent {\bf Proof.}} As in D. Zagier[Z] we introduce a real parameter with
respect to which we perform a peculiar induction to prove the following 
theorem, which in turn implies theorem (5.2) by the substitution 
 $\lambda = (d+\alpha(R))/n(R)$.
\vsp
{\bf Theorem \, 5.6}\, {\it 
Let the hypothesis be as in the previous theorem. The two
quantities 
 $$Q_{R,d}^{\lambda}(x) = \sum_{r \ge 1}\sum_{I} \sum_{\circ_{\lambda}} x^{
N_R(I;~d_{1},\ldots, d_{r})} \prod_{k=1}^r{Q_{R^I_{k},d_k}(x)}
 \eqno(5.7)$$
 $$ S_{R,d}^{\lambda}(x) = \sum_{r \ge 1}\sum_{I} 
\frac{x^{M'_R(I;~d)+M_R(I;~\lambda)}}{(x^{n(R^I_1)+n(R^I_2)}-1)
\ldots (x^{n(R^I_{r-1})+n(R^I_r)}-1)} \prod_{k=1}^r{Q_{R^I_{k}}(x)}  
\eqno(5.8)$$
 agree for every real number $\lambda \ge  (d+\alpha(R))/n(R)$.
\\ Here $\sum_{\circ_{\lambda}}$ denotes the summation over
$(d_{1},\ldots,d_{r}) \in \Z^r$ such that  $\sum_i{d_i}=d$ and
the following holds$$
\lambda \ge \frac {d_1 + \alpha(R^I_1)}{n(R^I_1)} > \frac {d_2 +
\alpha(R^I_2)}{n(R^I_2)} > \ldots > \frac {d_r + \alpha(R^I_r)}{n(R^I_r)}
\eqno(5.9)$$ 
}
\vsp
{\noindent {\bf Proof.}}
We first note that $Q_{R,d}^{\lambda}$ and 
$S_{R,d}^{\lambda}$ are step functions of $\lambda$ and they only 
jump at a discrete subset of $\R$. 
We assume by induction that $Q_{R',d}^{\lambda} =
S_{R',d}^{\la~}$ for all data $R'$ of rank $n(R') < n$, for all $d \in
\Z$
and $\lambda \in \R$. Now for a given data $R$ of rank $n(R) = n$ and
$d \in \Z$ we
make the following claims \\  
 {\bf Claim(1):} given any $N$, there exists $\lambda_0(N)$ such that
for $\lambda \ge \lambda_0(N)$, 
$Q_{R,d}^{\lambda}$ and $ S_{R,d}^{\lambda}$ 
agree modulo $x^N$. \\
{\bf Claim(2):} For any $\la \ge (d+\alpha(R))/n(R)$, if we define 
$Q_{R,d}^{\lambda^-}$ (resp. $S_{R,d}^{\lambda^-}$) to be 
$Q_{R,d}^{\lambda-\epsilon}$ (resp. $S_{R,d}^{\lambda-\epsilon}$)
 for $\epsilon > 0$ small enough such that the function
$Q_{R,d}^{\lambda}$ (resp. $S_{R,d}^{\lambda}$) has no jumps in the interval
$[\lambda-\epsilon, \lambda)$, then the two functions  
$\Delta Q_{R,d}^{\lambda}=
Q_{R,d}^{\lambda}-Q_{R,d}^{\lambda^-}$
 and 
$\Delta S_{R,d}^{\lambda}=S_{R,d}^{\lambda}-
S_{R,d}^{\lambda^-}$, are equal.
 \vs 
{\it Proof of claim(1):} For a particular $N$, by equation (5.1), 
the coefficient of $x^N$ in $Q_R$ involves
only finitely many choices of the integer $r$, partitions $I$, the
integers $(d_1,\ldots,d_r)$. Hence if we choose $\lambda_0(N) > (d_i +
\alpha(R^I_i))/n_i$ for all such combinations of $(r,I, d_1,\ldots,d_r)$,
 then the coefficient of $x^N$ in $Q_R$
 and $Q_{R,d}^{\lambda}$ are equal for $\lambda \ge \lambda_0(N)$. 
On the other hand, if $r>1$, we have $M_R(I;~\lambda)$
occurring in the exponent of the numerator which tend to $\infty$ as
$\lambda$ tends to $\infty$. So, for a fixed $N$ if we choose
$\lambda$ large enough, we do not get any contribution for the
coefficient of $x^N$ in $S_{R,d}^{\lambda}$. But for $r=1$ 
the part of the summation in  $S_{R,d}^{\lambda}$ is just $Q_{R}$.
hence the claim(1) follows.
\vs
{\it Proof of claim(2):}
\vs
Given a data $R$ and a sub-data $L$ of $R$
we define a numerical function 
$$\delta_R(L) = \sum_P \sum_{i>t}{(R-L)^P_iL^P_t}.  \eqno(5.10)$$

We first write down the recursions satisfied by the various numerical
functions that we have encountered in the statement of the theorem
(5.6).
\vsp
{\bf Lemma \, 5.11}\, {\it
Let $I$ be a partition of $R$ of length $r$. Let $0 < k < r$.
$$
\begin{array}{l}
{\bf (a)}~~\sigma_R(I) = \sigma_{R^I_{\le k}}(I_{\le k}) + 
\sigma_{R^I_{\ge k+1}}(I_{\ge k+1})
                                   + \delta_R(R^I_{\le k}) \\
 \\
{\bf (b)}~~N_R(I;~d_{1},\ldots, d_{r})- N_{R^I_{\ge 2}}(I_{\ge
                                    2};~d_2,\ldots ,d_r)  \\
~~~~~~~~~~~~~~~~~~~  = n(R^I_1)(n \la - d)-n \alpha(R^I_1)-
\delta_R(R^I_1) \\
 \\
{\bf (c)}~~M'_{R^I_{\le k}}(I_{\le k};~d(\la,R^I_{\le k}))+
 M'_{R^I_{\ge  k+1}}(I_{\ge k+1};~d-d(\la,R^I_{\le k})) \\ 
    \\
 ~~~~~~~~~~=M'_R(I;~d)-(2n(R^I_{\le k})-n(R)-n(R^I_k)+n(R^I_r))d(\la,R^I_{\le k}) \\
 ~~~~~~~~~~~~~~~~~~~~~~~~~~~ - n(R^I_k)- n(R^I_{k+1})+\delta_R(R^I_{\le
k}) + n(R^I_{\le k})d.
\end{array}
\eqno(5.12)$$
where 
$d(\la,L)= n(L) \la - \alpha(L)$ for any data $L$.
}
\vsp
{\noindent {\bf Proof.}} All these statements follow from  straight
forward calculations, so we will not give the details.  $\hfill{\Box}$
\vsp
 We now compute $\Delta Q_{R,d}^{\lambda}$. It is zero unless
there is a partition $I$ of $R$ and a $r$-tuple $(d_1,\ldots,d_r)$ with
$\sum d_i = d$ such that  $\la = (d_1 + \alpha(R^I_1))/n(R^I_1)$. 
 For such a $\la$, we observe that 
$$\Delta Q_{R,d}^{\lambda} = \sum_{r \ge 1}\sum_{I}\sum_{\stackrel
{\circ_{\lambda}}{ \la = (d_1 + \alpha(R^I_1))/n(R^I_1)}}
 x^{N_R(I;~d_{1},\ldots, d_{r})} \prod_{k=1}^r{Q_{R^I_k,d_k}(x)} \eqno(5.13)
$$
We can use the lemma (5.11)  in the above formula  and separate the
expressions which have $k=1$ and $k \ge 2$. Hence the right hand side in the
equation (5.13) becomes

$$
\sum_{\stackrel{L~\mbox{{\footnotesize sub-data of}}~R}{d(\la,L) \in \sZ}}
 x^{ n(L)(n(R) \la - d)-n(R) \alpha(L)- 
\delta_R(L)} Q^{\lambda}_{L, d(\la,L)}
Q_{R-L,d - d(\la,L)}^{\la^-}. \eqno(5.14)  
$$
 
  Now we compute $\Delta S_{R,d}^{\lambda}$ at a $\la$ when there
is a jump. This happens when 
$(n(R^I_1)+\ldots +n(R^I_k))\la - \alpha(R^I_k)$ is an integer for some
partition $I$ of $R$  with $\ell (I)=r$ and for some positive integer $k < r$. 
\vs
Fix a partition $I$ of length $r$.
Let $$ \pi_I = \{k<r| d(\la, R^I_{\le k}) \in \Z \} \eqno(5.15)$$
 One can see that 
  $\Delta M(I;~\la) = \sum_{k \in \pi_I} (n(R^I_k)+n(R^I_{k+1}))$ so 
$$
\begin{array}{l} 
 x^{M(I;~\la)} -x^{M(I;~\la^-)}=x^{M(I;~\la^-)}(x^{\sum_{k \in \pi_I}
 (n(R^I_k)+n(R^I_{k+1}))}-1) \\
 \\
 ~~~~~~~~= \sum_{k \in \pi_I}x^{M(I;~\la^-)+ \sum_{\{k' \in \pi_I|k' < k\}}
(n(R^I_{k'})+n(R^I_{k'+1}))} (x^{(n(R^I_k)+n(R^I_{k+1}))}-1)\\
 \\
~~~~~~~~= \sum_{k \in \pi_I}x^{M(R^I_{\le k};~\la)+M(R^I_{\ge k+1};~\la^-)+
(2n(R)-2n(R^I_{\le k})+n(R^I_k)-n(R^I_r))d(\la,R^I_{\le k})} \\
~~~~~~~~~~~~~~~~~~~~~~~~~~~~~~~.(x^{(n(R^I_k)+n(R^I_{k+1}))}-1)
\end{array}
\eqno(5.16)$$
Let $g_R(I;~d)$ denote the following rational function of $x$ 
     $$\frac{x^{M'_R(I;~d)+M_R(I;~\la)}}{(x^{n(R^I_1)+n(R^I_2)}-1)
\ldots (x^{n(R^I_{r-1})+n(R^I_r)}-1)} \eqno(5.17)$$
 Using the lemma (5.11) and equation (5.16) we can verify that 
$\Delta(g_R(I;~d))$ \\is equal to the following
$$ 
\begin{array}{l}
\sum_{k \in \pi_I}x^{ n(R^I_k)(n(R) \la - d) 
-n(R) \alpha(R^I_{\le k})- \delta_R(R^I_{\le k})}\\
~~~~~~~~~~~~.g_{R^I_{\le k}}(I_{\le k},d(\la,R^I_{\le k}))
g_{R^I_{\ge k+1}}(I_{\ge k+1},d-d(\la,R^I_{\le k}))
\end{array}
 \eqno(5.18)$$

Now $\Delta S_{R,d}^{\lambda}$ is computed to be
$$
\sum_{r \ge 1}\sum_{I} 
\Delta(g_R(I;~d)) \prod_{k=1}^r{Q_{R^I_k}(x)} 
\eqno(5.19)$$
Using the equation (5.18), and grouping together all terms which give the 
sub-data $L$, we get the following expression for 
$\Delta S_{R,d}^{\lambda}$
$$ 
\sum_L
 x^{ n(L)(n(R) \la - d)-n(R) \alpha(L)- \delta_R(L)}
S_{L, d_{\la,L}}^{\la}S_{R-L,d - d(\la,L)}^{\la^-}. \eqno(5.20)  
$$     
 where the summation is over sub-data $L$ of $R$ with $ d(\la,L) \in \Z$. 
But \\$Q_{L, d_{\la,L}}^{\la} = S_{L, d_{\la,L}}^{\la}$ and 
$Q_{R-L,d - d(\la,L)}^{\la^-}=S_{R-L,d - d(\la,L)}^{\la^-}$
by induction (since $n(L)$ and $n(R)-n(L)$ are less than $n(R)$), hence we
get $\Delta S_{R,d}^{\lambda}=\Delta Q_{R,d}^{\lambda}$.
this proves claim(2).
 \vs      
To prove the theorem it is enough to check that the coefficient of
$x^N$ in $Q_{R,d}^{\lambda}$ and in $S_{R,d}^{\lambda}$ agree
for any $N$.
  For a given  $N$, the claim(1) implies that 
the coefficients of $Q_{R,d}^{\lambda}$ and $S_{R,d}^{\lambda}$ are
equal when $\la$ is sufficiently large . Since $Q_{R,d}^{\lambda}$ and 
$S_{R,d}^{\lambda}$  are step functions
of $\lambda$ jumping only at a discrete set of real numbers, and for
such real numbers by claim(2) their jumps agree therefore the jumps in
the coefficients also agree, which in turn proves that the coefficients
are the same. This completes the proof of the theorem.   $\hfill{\Box}$
\vsp
Now if we define $$\sigma'_R(I)= \sum_{P\in S} 
    \sum_{k>l,i<t}{I^P_{i,k}I^P_{t,l}} \eqno(5.21)$$ 
then  one observes that dimensions of the flag varieties ${\cal F}_R$ and 
${\cal F}_{R^I_k}$ are related by 
$$ \dim {\cal F}_R - \sum_{k=1}^r\dim {\cal F}_{R^I_k}=
\sigma_R(I)+\sigma'_R(I). \eqno(5.22)$$ 
Using this expression we can formulate the closed formula for the 
Poincar\'e polynomial of the moduli space of \pssb s as follows.
\vsp
{\bf Theorem\, 5.23.}\, {\it
 The Poincar\'e polynomial $P_{R,d}$ of the moduli space
of parabolic stable bundles with a fixed
determinant of degree $d$, and data $R$ satisfying the condition `par
semi-stable = par stable' is given by 

$$\frac{1-t^2}{(1+t)^{2g}}\sum_{r \ge 1}\sum_{I} 
\frac{t^{2(\sigma'_R(I)-(n(R)-n(R^I_r))d+M_g(I;~(d+\alpha(R))/n(R))}}
{(t^{2n(R^I_1)+2n(R^I_2)}-1)\ldots (t^{2n(R^I_{r-1})+2n(R^I_r)}-1)} 
\prod_{k=1}^r{P_{R^I_k}(t)} \eqno(5.24)$$
 where $M_g(I;~\lambda)$ is
$$
\begin{array}{r}
\sum_{k=1}^{r-1}(n(R^I_k)+n(R^I_{k+1}))
([(n(R^I_1)+\ldots + n(R^I_k))\lambda - \alpha(R^I_{\le k})]+1) \\
+(g-1)\sum_{i<j}n(R^I_i)n(R^I_j)
\end{array}
\eqno(5.25)$$
and $P_{R}(t)$ is defined to be 
 $$\left( \frac{\prod_{i=1}^{n(R)}{(1-t^{2i})}^{|S|}}
{\prod_{P\in S}\prod_{\{i| {R^P_i} \ne
0\}}\prod_{l=1}^{R^P_i}{(1-t^{2l})}} \right)
\left( \frac{\prod_{i=1}^{n(R)}{(1+t^{2i-1})^{2g}}}
{(1-t^{2n(R)})\prod_{i=1}^{n(R)-1}{(1-t^{2i})^2}}\right).\eqno(5.26)$$
}
\section{Sample calculations}
\medskip
\centerline{\bf Rank 2 }
\smallskip 
Now we write down the Poincar\'e polynomial in more and more explicit forms
for any data $R$ such that $n(R)=2$.
\vs
Let $T$ be a subset of $S$ defined by $\{P \in S|R^P_1 = 1\}$, which
is the set of vertices where the parabolic filtration is non-trivial. 
Then we get $$P_R(t) = \frac {(1+t^2)^{|T|}(1+t)^{2g}
(1+t^3)^{2g}}{(1-t^4){(1-t^2)}}  \eqno(6.1)$$

Given any partition $I$ of $R$, we define a subset $T_I$ of $T$ by
$$T_I=\{P \in T|I^P_{1,1} =0\}  \eqno(6.2)$$
from this definition we observe that ${\sigma}'_R(I)$ (as defined in 
 (5.21)) is just $|T_I|$.
 \vs 
Let $\chi_I : T~ \longrightarrow ~\{ 1,-1\}$ be defined by
$$ 
\begin{array}{lll}
 \chi_I(P) & =~~1 &  \mbox{if}~~ P \in T_I \\
           & =-1 &  \mbox{otherwise}
\end{array}
\eqno(6.3)$$
Using the theorem (5.23) for rank 2 moduli we obtain the following.
\vsp
 {\bf Proposition\, 6.4.}\, {\it
For any degree $d$, the Poincar\'e polynomial for the moduli
space of rank 2 parabolic bundles  with data $R$ and satisfying the 
condition `par semi-stable = par stable' is given by 
$$ P_{R,d}(t)= \frac
{(1+t^2)^{|T|}(1+t^3)^{2g}-(\sum_It^{2(g+|T_I|+[\psi_I]+a_I)
})(1+t)^{2g}}{(1-t^4){(1-t^2)}} \eqno(6.5)$$
 where 
$$\psi_I = \sum_{P \in T}\chi_I(P)(\alpha^P_1 -\alpha^P_2) \eqno(6.6)$$
and  $a_I$ is $1$ or $0$ depending on whether $d+ [\psi_I]$ is even or
odd.
}
\vsp
Now we put $g=0$ in the formula. 
Since $P_{R,d}$ is a power series in $t$, one sees that 
$|T_I|+[\psi_I]+a_I \ge 0$ for every partition $I$    
\vs
Using the proposition (6.4), 
the zeroth Betti number of the moduli space can be computed to be
equal to $1 - |\{I||T_I|+[\psi_I]+a_I =0\}|$, hence the quantity 
$|T_I|+[\psi_I]+a_I$ is 0 for at most one partition. 
Hence we obtain the following corollary 
\vsp
{\bf Corollary \, 6.7.}\, {\it Assuming the condition \ps= ~we have \\
 {\bf a)} The moduli space of \pssb s of rank 2 is non-empty iff for 
every partition $I$ we have $|T_I|+[\psi_I]+a_I >0$.\\
{\bf b)} The moduli is actually connected when it is non-empty. 
}
\vsp
One can easily see that this condition is equivalent to the condition
given by I.Biswas [B].
Even in higher rank we can get a criterion for existence of stable
bundles by setting $P_{R,d}(t) \ne 0$. 
\vs
In what follows we assume that $\psi_I$ is never an integer, which has
the effect that the condition \ps= ~holds for all the degrees.
\vs
Using the above formula  for the Poincar\'e polynomial we compute
it in a explicit form, when the cardinality of $S$ is small (1,2,3 and 4).
For this, one observes that the above expression for $ P_{R,d}(t)$, the
dependence on the weights is only via their differences. In view of this
we define $\delta^P = \alpha^P_1 -\alpha^P_2$ for each $P \in S$.
\vs
When $S=\{P \}$, $\delta^P$ arbitrary, $R^P_i =1$ for all $i$ and any degree
$d$, we compute the Poincar\'e polynomial to be
 $$ P_{R,d}(t) = \frac{(1+t^3)^{2g}-t^{2g}(1+t)^{2g}}{(1-t^2)^2}. \eqno(6.8) $$

When $S=\{ P_1,P_2\}$, $\delta^{P_1}$ and $\delta^{P_2}$ arbitrary, 
$R^P_i =1$ for all $i$ and $P$, and any degree
$d$, we have
 $$ P_{R,d}(t) =
\frac{(1+t^2)((1+t^3)^{2g}-t^{2g}(1+t)^{2g})}{(1-t^2)^2}.
\eqno(6.9) $$
 
When $S= \{P_1,P_2,P_3 \}$, $R^{P_j}_i=1$ for all $i$ and $j= 1,\ldots
 3$, and any degree $d$. By reordering $P_1,P_2,P_3$, we may assume that 
$\delta^{P_1} \le \delta^{P_2} \le \delta^{P_3} $. Now there are two 
possibilities \\
(i) If $\delta^{P_1}+ \delta^{P_2}+ \delta^{P_3} <-2$ or  $-\delta^{P_1}+ \delta^{P_2}+ \delta^{P_3} >0$
then 
 $$ P_{R,d}(t)
=\frac{(1+t^2)^2((1+t^3)^{2g}-t^{2g}(1+t)^{2g})}{(1-t^2)^2}. 
\eqno(6.10)
$$
(ii) If $\delta^{P_1}+ \delta^{P_2}+ \delta^{P_3} >-2$ and
$-\delta^{P_1}+ \delta^{P_2}+ \delta^{P_3} <0$ (which is the remaining
case),  then 
$$
P_{R,d}(t)=\frac{(1+t^2)^2(1+t^3)^{2g}-4t^{2g+2}(1+t)^{2g}}{(1-t^2)^2}.
\eqno(6.11)
$$
 
When $S= \{P_1,P_2,P_3,P_4 \}$, $R^{P_j}_i =1$ for all $i$ and $j$, 
and any degree $d$. 
Again we assume $\delta^{P_1} \le \delta^{P_2} \le \delta^P_3 \le
\delta^{P_4}$.
Again there are two possibilities
\vs
(i) If $\delta^{P_1}+ \delta^{P_2}+ \delta^{P_3} -\delta^{P_4}<-2$ or  $-\delta^{P_1}+ \delta^{P_2}+ \delta^{P_3}+\delta^{P_4} >0$
then 
 $$ P_{R,d}(t)
=\frac{(1+t^2)^3((1+t^3)^{2g}-t^{2g}(1+t)^{2g})}{(1-t^2)^2}.   
\eqno(6.12)
$$
(ii) If $\delta^{P_1}+ \delta^{P_2}+ \delta^{P_3} -\delta^{P_4}>-2$ and  
 $-\delta^{P_1}+ \delta^{P_2}+ \delta^{P_3}+\delta^{P_4} <0$ (which is
the remaining case),  then 
$$
P_{R,d}(t)=\frac{(1+t^2)^3(1+t^3)^{2g}-4t^{2g+2}(1+t^2)(1+t)^{2g}}{(1-t^2)^2}.
\eqno(6.13)
$$
\vsp
{\bf Remark\, 6.14.}\, 
 Note that in these cases we have considered, the
Poincar\'e polynomial does not depend on the degree. In fact we can
verify that in general for rank 2 the Poincar\'e polynomial is
independent of the degree.
\vsp
\centerline{\bf Rank 3 and 4 }
\medskip
Using the software Mathematica, we have computed the Poincar\'e
polynomials and the Betti numbers for the rank 3 and rank 4 when the
number of parabolic points is one or two.
In this case we find that  the Poincar\'e polynomial has
dependence on the weights and degree. In the appendix we actually give the
tables for the Betti numbers (in the rank 3 and rank 4 case) taking 
different set of weights into consideration. 
For the Poincar\'e polynomial we choose one set of
weights as an example in each of the following cases.

 \vs
When rank $=3$, $S=\{P_1\}$,  $R^{P_1}_i=1$ for all $i$ and assume 
that the condition \ps= ~holds. Then for 
all choices of weights and degree $d$ we have 
$$ 
\begin{array}{l}
P_{R,d}= \{ t^{6g-2}(1+t^{2}+t^{4})\left( 1 + t \right)^{4g}
 -t^{4g-2}(1+t^2)^2 \left( 1 + t \right)^{2g} 
      \left( 1 + {t^3} \right)^{2g} \\
~~~~~~~~~~ +\left( 1 + {t^3} \right)^{2g}
      \left( 1 + {t^5} \right)^{2g}\}/(       
     \left({t^2}-1 \right)^4 \left( 1 + {t^2} \right)).
\end{array}
\eqno(6.15)
$$
When rank $=3$, $S=\{P_1,\,P_2 \}$,  $R^{P_j}_i=1$ for all $i$ and
$j$. One observes that the condition \ps= ~ holds for all choices of 
degree.
 For $(\alpha^{P_1}_1,\,\alpha^{P_1}_2,\,\alpha^{P_1}_3)=(0,\,1/12 ,\,3/12 )$,
$(\alpha^{P_2}_1,\,\alpha^{P_2}_2,\,\alpha^{P_2}_3)=(1/12,\,5/12,\,6/12)$, one
observes that the condition \ps= ~ holds for all choices of 
degree.
When the degree $d=0~\mbox{or}~2$ mod 3 we find that
$$ 
\begin{array}{l}
P_{R,d}=\{ -3t^{4g}\left( 1 + t \right)^{2g}
        \left( 1 + {t^2} \right)^2
       \left( 1 + {t^3} \right)^{2g} +t^{6g}\left( 1 + t \right)^{4g}
        \left( 2 + 5{t^2} + 2{t^4} \right) \\ 
~~~~~~~~~~ + \left( 1 + {t^3} \right)^{2g}
       \left( 1 + {t^5} \right)^{2g}\left( 1 + t^2 +t^4 \right)\}/
       \left( 1 - {t^2} \right)^4
\end{array}
\eqno(6.16)
$$
and if $d=1$ mod 3 then the Poincar\'e polynomial is
$$ 
\begin{array}{l}  
P_{R,d}= \left( 1 + {t^2} + {t^4} \right) 
     \{ t^{6g-2}\left( 1 + {t^2} + {t^4} \right)
                 \left( 1 + t \right)^{4g}\\ 
~~~~~~- t^{4g-2}\left( 1 + {t^2} \right)^2 \left( 1 + t \right)^{2g}
        \left( 1 + {t^3} \right)^{2g}+ \left( 1 + {t^3} \right)^{2g}
        \left( 1 + {t^5} \right)^{2g} \}/\left(1- {t^2} \right)^4.
\end{array}
\eqno(6.17)
$$
When rank is 4 and $|S|=1$ we find that the Poincare Polynomial depends on the
degree too. If we choose $R^P_i=1$ for all $i$, and choose 
$(\alpha^P_1,\,\alpha^P_2,\,\alpha^P_3,\,\alpha^P_4)=(0,\,1/8,\,1/4,\,1/2)$,
then the condition \ps= ~ holds for all choices of degree.
We find that that 
$$ 
\begin{array}{l}
P_{R,0}=P_{R,1}=P_{R,2}=\{ {\left( 1 + {t^3} \right)^{2g}
       \left( 1 + {t^5} \right)^{2g}
       \left( 1 + {t^7} \right)^{2g}} \\ 
~~~~~~- {2t^{-2 + 6g}\left( 1 + t \right)^{2g}
       \left( 1 + {t^3} \right)^{2g} \left( 1 + {t^5} \right)^{2 g}
       \left( 1 + {t^2} + {t^4} \right)}  \\ 
~~~~~~ - {t^{-4 + 8g}\left( 1 + t \right)^{2g}
       \left( 1 + {t^3} \right)^{4g}\left( 1 + {t^2} + {t^4}
             \right)^2} \\  
~~~~~~ + {t^{-4 + 10 g}\left( 1 + {t^2} \right) \left( 1 + t \right)^{4 g} 
     \left( 1 + {t^3} \right)^{2 g} 
     \left( 3 + 5 {t^2} + 5 {t^4} + 3 {t^6} \right)}  \\ 
~~~~~~ -{2 t^{-4 + 12 g} \left( 1 + t \right)^{6 g} 
     \left( 1 + {t^2} + {t^4} \right)^2 } \}/(\left( 1 - {t^2} \right)^6
      \left( 1 + {t^2} \right) \left( 1 + {t^2} + {t^4} \right))
\end{array}
\eqno(6.18)
$$ 
and
$$
\begin{array}{l}
P_{R,3}= \{ {\left( 1 + {t^3} \right)^{2g}
       \left( 1 + {t^5} \right)^{2g}
       \left( 1 + {t^7} \right)^{2g}} \\ 
~~~~ -{t^{-4 + 6g} \left( 1 + t \right)^{2g}
       \left( 1 + {t^3} \right)^{2g} \left( 1 + {t^5} \right)^{2 g}
       \left( 1 + {t^2} + {t^4} \right) \left(1+t^4 \right)} \\ 
~~~~ -{t^{-4 + 8g}\left( 1 + t \right)^{2g}
       \left( 1 + {t^3} \right)^{4g}\left( 1 + {t^2} + {t^4} \right)^2}\\ 
~~~~ +{t^{-6 + 10 g}\left( 1 + {t^2} \right)^4 \left( 1 + t \right)^{4 g} 
        \left( 1 + {t^3} \right)^{2 g} \left( 1 + {t^4} \right)}  \\ 
 -{2t^{-6 + 12 g} \left( 1 + t \right)^{6 g} \left( 1 + {t^4} \right) 
     \left( 1 + {t^2} +{t^4} \right)^2 } \}/({\left( 1 - {t^2} \right)^6
      \left( 1 + {t^2} \right) \left( 1 + {t^2} + {t^4} \right)}).
\end{array}
\eqno(6.19)
$$
If we choose the weights $(\alpha^P_1,\,\alpha^P_2,\,
\alpha^P_3,\,\alpha^P_4)=(0,\,1/5,\,4/5,\,9/10)$ then again the condition \ps=
~ holds for all choices of degree. If $P_{R,d}$ and $P'_{R,d}$ denote
the Poincar\'e polynomial for the moduli space of parabolic stable
bundles with with data $R$ (satisfying n(R)=4), having degree $d$
and with  weights (0,\,1/8,\,1/4,\,1/2) and (0,\,1/5,\,4/5,\,9/10) respectively, then
we find that
$ P'_{R,0}=P'_{R,2}= P_{R,0}~~\mbox{and}~~P'_{R,1}=P'_{R,3}= P_{R,1}$.
\newpage
\section{Appendix : Betti number tables}
The following tables give the Betti numbers up to the middle dimension
of the moduli space of parabolic bundles over $X$ for rank 2, 3 and 4
and low genus. When $\beta_0 =0$, we mean that the space is empty. 

\centerline{\bf Rank 2}
\medskip
Any degree $d$, $R^P_i =1$ for all $i$ and $P \in S$.

{\bf Case A)} $S=\{P_1\}$, $\delta^{P_1}$ arbitrary.\\
{\bf Case B)} $S= \{P_1,\,P_2\}$, $\delta^{P_1}$ and $\delta^{P_2}$ arbitrary.\\
{\bf Case C)} $S= \{P_1,\,P_2,\,P_3 \}$, $\delta^{P_1}+ \delta^{P_2}+ 
\delta^{P_3} <-2$ or \\ 
$-\delta^{P_1}+ \delta^{P_2}+ \delta^{P_3} >0$\\
{\bf Case D)} $S= \{P_1,\,P_2,\,P_3 \}$, $\delta^{P_1}+ \delta^{P_2}+ 
\delta^{P_3} >-2$ and \\ 
$-\delta^{P_1}+ \delta^{P_2}+ \delta^{P_3} <0$.\\
{\bf Case E)} $S= \{P_1,\,P_2,\,P_3,\,P_4 \}$, $\delta^{P_1}+ \delta^{P_2}+ 
\delta^{P_3} -\delta^{P_4}<-2$ or \\ 
$-\delta^{P_1}+ \delta^{P_2}+ \delta^{P_3}+\delta^{P_4} >0$\\
{\bf Case F)} $S= \{P_1,\,P_2,\,P_3,\,P_4 \}$,\, $\delta^{P_1}+ \delta^{P_2}+
 \delta^{P_3} -\delta^{P_4}>-2$ and \\ 
$-\delta^{P_1}+ \delta^{P_2}+ \delta^{P_3}+\delta^{P_4} <0$.
{\small
\begin{center}
\begin{tabular}{|l||l|l|l|l|l|l||l|l|l|l|l|r||}  \hline 
 & \multicolumn{6}{c||}{Genus g=0} &\multicolumn{6}{c||}{Genus g=1}
\\ \cline{2-13}
  &A&B&C&D&E&F&A&B&C&D&E&F  \\ \hline
$\beta_0   $&0&0&0~&1~&0~&1 &1 &1 &1 &1 &1  &1     \\ \hline
$\beta_1   $&-&-&- &- &- &0 &0 &0 &0 &0 &0  &0     \\ \hline
$\beta_2   $&-&-&- &- &- &- &- &2 &3 &4 &4  &5     \\ \hline
$\beta_3   $&-&-&- &- &- &- &- &- &0 &2 &0  &2     \\ \hline 
$\beta_4   $&-&-&- &- &- &- &- &- &- &- &6  &8     \\ \hline \hline
 & \multicolumn{6}{c||}{Genus g=2} &\multicolumn{6}{c||}{Genus g=3}
\\ \cline{2-13}
 &A&B&C&D&E&F&A&B&C&D&E&F  \\ \hline
$\beta_0   $&1&1&1 &1 &1 &1 &1 &1 &1 &1 &1  &1      \\ \hline
$\beta_1   $&0&0&0 &0 &0 &0 &0 &0 &0 &0 &0  &0      \\ \hline
$\beta_2   $&2&3&4 &4 &5 &5 &2 &3 &4 &4 &5  &5      \\ \hline
$\beta_3   $&4&4&4 &4 &4 &4 &6 &6 &6 &6 &6  &6      \\ \hline 
$\beta_4   $&2&4&7 &8 &11&12&3 &5 &8 &8 &12 &12     \\ \hline
$\beta_5   $&-&8&12&16&16&20&12&18&24&24&30 &30     \\ \hline
$\beta_6   $&-&-&8 &14&15&22&18&21&26&27&34 &35     \\ \hline
$\beta_7   $&-&-&- &- &24&32&12&24&42&48&66 &72     \\ \hline
$\beta_8   $&-&-&- &- &- &- &- &36&57&72&83 &99     \\ \hline
$\beta_9   $&-&-&- &- &- &- &- &- &48&68&90 &116    \\ \hline
$\beta_{10}$&-&-&- &- &- &- &- &- &- &- &114&144    \\ \hline
\end{tabular}
\end{center}
}
\newpage
\centerline{\bf Rank 3}
\medskip
{\bf Case A)}$S=\{ P\}$, $R^P_i=1$ for all $i$.\\ 
We take all choices of weights and degrees.\\
{\bf Case B)} When $S=\{ P_1,\,P_2 \}$, $R^{P_1}_i=1=R^{P_2}_i$ 
for all $i$.\\   
$(\alpha^{P_1}_1,\,\alpha^{P_1}_2,\,\alpha^{P_1}_3)$= (0,\,1/12,\,3/12),
$(\alpha^{P_2}_1,\,\alpha^{P_2}_2,\,\alpha^{P_2}_3)$= (1/12,\,5/12,\,6/12)\\
d=0 or 2 mod 3\\
{\bf Case C)} $S=\{ P_1,\,P_2 \}$, $R^{P_1}_i=1=R^{P_2}_i$ for all $i$.\\ 
$(\alpha^{P_1}_1,\,\alpha^{P_1}_2,\,\alpha^{P_1}_3)$= (0,\,1/12,\,3/12), 
$(\alpha^{P_2}_1,\,\alpha^{P_2}_2,\,\alpha^{P_2}_3)$= (1/12,\,5/12,\,6/12),\\
 d=1 mod 3.\\
[2mm]
\begin{tabular}{|l||l|l|l|l||l|l|l|l||l|l|l|r||}  \hline  
 & \multicolumn{4}{c||}{A,\,g\,=} &\multicolumn{4}{c||}{B,\,g\,=}&\multicolumn{4}{c||}{C,\,g\,=}     \\ \cline{2-13}
            &0 &1&2 &3    &0&1 &2  &3     &0&1 &2  &3     \\ \hline \hline
$\beta_0   $&0 &1&1 &1    &0&1 &1  &1     &0&1 &1  &1      \\ \hline
$\beta_1   $&- &0&0 &0    &-&0 &0  &0     &-&0 &0  &0       \\ \hline
$\beta_2   $&- &2&3 &3    &-&5 &5  &5     &-&4 &5  &5       \\ \hline
$\beta_3   $&- &0&4 &6    &-&2 &4  &6     &-&0 &4  &6       \\ \hline
$\beta_4   $&- &-&7 &7    &-&12&15 &15    &-&8 &15 &15      \\ \hline
$\beta_5   $&- &-&16&24   &-&6 &24 &36    &-&0 &24 &36      \\ \hline
$\beta_6   $&- &-&18&28   &-&16&40 &49    &-&10&39 &49      \\ \hline
$\beta_7   $&- &-&36&60   &-&- &80 &120   &-&- &76 &120     \\ \hline
$\beta_8   $&- &-&45&103  &-&- &108&176   &-&- &98 &176     \\ \hline
$\beta_9   $&- &-&56&140  &-&- &188&314   &-&- &164&314     \\ \hline
$\beta_{10}$&- &-&70&261  &-&- &251&531   &-&- &203&530     \\ \hline
$\beta_{11}$&- &-&64&354  &-&- &344&784   &-&- &264&778     \\ \hline
$\beta_{12}$&- &-&- &537  &-&- &436&1312  &-&- &318&1293    \\ \hline
$\beta_{13}$&- &-&- &780  &-&- &480&1878  &-&- &332&1828    \\ \hline
$\beta_{14}$&- &-&- &998  &-&- &528&2816  &-&- &370&2697    \\ \hline
$\beta_{15}$&- &-&- &1380 &-&- &-  &4036  &-&- &-  &3788    \\ \hline
$\beta_{16}$&- &-&- &1652 &-&- &-  &5454  &-&- &-  &4983    \\ \hline
$\beta_{17}$&- &-&- &1936 &-&- &-  &7442  &-&- &-  &6610    \\ \hline
$\beta_{18}$&- &-&- &2170 &-&- &-  &9346  &-&- &-  &8007    \\ \hline
$\beta_{19}$&- &-&- &2160 &-&- &-  &11526 &-&- &-  &9572    \\ \hline
$\beta_{20}$&- &-&- &-    &-&- &-  &13394 &-&- &-  &10812   \\ \hline
$\beta_{21}$&- &-&- &-    &-&- &-  &14562 &-&- &-  &11508   \\ \hline
$\beta_{22}$&- &-&- &-    &-&- &-  &15210 &-&- &-  &11984   \\ \hline
\end{tabular}
\newpage
\centerline{\bf Rank 4}
\medskip
$|S|=1$ \\
{\bf Case A)}  $R^P_i=1$ for all $i$, d= 0 or 1 or 2 mod 4,\\
$(\alpha^P_1,\alpha^P_2,\alpha^P_3,\alpha^P_4)=(0,\,1/8,\,1/4,\,1/2)
$~~~~~ Or\\
$R^P_i=1$ for all $i$, d= 0 or 2 mod 4,\\
$(\alpha^P_1,\,\alpha^P_2,\,\alpha^P_3,\,\alpha^P_4)=(0,\,1/5,\,4/5,\,9/10)$.\\
{\bf Case B)}  $R^P_i=1$ for all $i$, d= 3 mod 4,\\ 
$(\alpha^P_1,\,\alpha^P_2,\,\alpha^P_3,\,\alpha^P_4)=
(0,\,1/8,\,1/4,\,1/2)$~~~~~ Or\\
$R^P_i=1$ for all $i$, d= 1 or 3 mod 4,\\ 
$(\alpha^P_1,\,\alpha^P_2,\,\alpha^P_3,\,\alpha^P_4)=
(0,\,1/5,\,4/5,\,9/10)$. \\
[2mm]
\begin{tabular}{|l||l|l|l||l|l|r||}   \hline  
 & \multicolumn{3}{c||}{A,g} &\multicolumn{3}{c||}{B,g}  \\ \cline{2-7}
            &0&1 &2    &0&1  &2    \\ \hline \hline
$\beta_0   $&0&1 &1    &0&1  &1    \\ \hline
$\beta_1   $&-&0 &0    &-&0  &0    \\ \hline
$\beta_2   $&-&4 &4    &-&3  &4    \\ \hline
$\beta_3   $&-&2 &4    &-&0  &4    \\ \hline
$\beta_4   $&-&8 &11   &-&5  &11   \\ \hline
$\beta_5   $&-&4 &20   &-&0  &20   \\ \hline
$\beta_6   $&-&10&31   &-&6  &31   \\ \hline
$\beta_7   $&-&- &64   &-&-  &64   \\ \hline
$\beta_8   $&-&- &90   &-&-  &89   \\ \hline
$\beta_9   $&-&- &164  &-&-  &160  \\ \hline
$\beta_{10}$&-&- &241  &-&-  &232  \\ \hline
$\beta_{11}$&-&- &376  &-&-  &356  \\ \hline
$\beta_{12}$&-&- &563  &-&-  &521  \\ \hline
$\beta_{13}$&-&- &792  &-&-  &712  \\ \hline
$\beta_{14}$&-&- &1144 &-&-  &1001 \\ \hline
$\beta_{15}$&-&- &1508 &-&-  &1272 \\ \hline
$\beta_{16}$&-&- &2003 &-&-  &1635 \\ \hline
$\beta_{17}$&-&- &2492 &-&-  &1952 \\ \hline
$\beta_{18}$&-&- &2989 &-&-  &2263 \\ \hline
$\beta_{19}$&-&- &3424 &-&-  &2528 \\ \hline
$\beta_{20}$&-&- &3675 &-&-  &2660 \\ \hline
$\beta_{21}$&-&- &3816 &-&-  &2760 \\ \hline
$\beta_{22}$&-&- &-    &-&-  &-    \\ \hline
\end{tabular}

\section*{References} 
[A-B] Atiyah, M. F. and Bott, R. : The Yang-Mills equations 
over Riemann surfaces.
{\sl Philos. Trans. Roy. Soc. London} Series A, {\bf 308} (1982) 523-615

[B] I. Biswas : A criterion for the existence of parabolic stable
bundle of rank two over the projective line, to appear in the {\sl 
Int. Jour. math.} {\bf 9} (1998), 523-533. 

[D-R] Desale, U. V. and Ramanan, S. : Poincar\'e Polynomials of
the Variety of Stable Bundles, {\sl Math. Annln.} {\bf 216}
(1975), 233-244.

[F-S] Furuta, M. and Steer, B. : Siefert-fibered homology
3-spheres and Yang-Mills equations on Riemann surfaces with
marked points, {\sl Adv. Math.} {\bf 96} (1992) 38-102. 

[G-L] Ghione, F. and Letizia, M. : Effective divisors of higher
rank on a curve and the Siegel formula, {\sl Composito Math.}
{\bf 83} (1992), 147-159.

[H-N] Harder, G. and Narasimhan, M. S. : On the Cohomology
Groups of Moduli Spaces of Vector Bundles over Curves, {\sl
Math. Annln.} {\bf 212} (1975), 215-248.

[M-S] Mehta, V. B. and Seshadri, C. S. : Moduli of vector bundles
on curves with parabolic structures, {\sl Math. Annln.} {\bf 248}
(1980) 205-239.

[N1] Nitsure, N. : Cohomology of the moduli of parabolic vector
bundles, {\sl Proc. Indian Acad. Sci. (Math. Sci.)} {\bf 95}
(1986) 61-77.

[N2] Nitsure, N. : Quasi-parabolic Siegel formula. Proc. Indian Acad. Sci.
(Math. Sci.) {\bf 106} (1996) 133-137, Erratum: {\bf 107} (1997) 221-222. 
(alg-geom/9503001 on the Duke e-print server.)

[S] Seshadri, C. S. : Fibres vectoriels sur les courbes
algebriques, {\sl Asterisque} {\bf 96} (1982).

[Z] Zagier, Don. : Elementary aspects of Verlinde formula and of the 
Harder-Narasimhan-Atiyah-Bott formula, {\sl Israel
Math. Conf. Proc.} {\bf 9} (1996) 445-462.   

\bigskip
Address: \begin{flushleft}
Tata Institute of Fundamental Research,\\
Mumbai 400 005,\\
India.\\
\end{flushleft}

\end{document}